\documentclass[12pt,a4paper]{amsart}

\usepackage[utf8]{inputenc}
\usepackage[T1]{fontenc}
\usepackage{lmodern}
\usepackage{microtype}
\usepackage[svgnames]{xcolor}
\usepackage{tikz}
\usetikzlibrary{arrows,fadings,patterns}

\usepackage {amssymb}
\usepackage {amsmath}
\usepackage {amsthm}
\usepackage{mathtools}
\usepackage{xfrac}
\usepackage{cancel}

\usepackage{url}

\usepackage{enumitem}
\usepackage{hyperref}
\hypersetup{
    colorlinks=true,
    linkcolor=blue,
    filecolor=magenta,
    urlcolor=cyan,
    citecolor=magenta
}
\usepackage{cleveref}

\DeclareMathOperator {\ord} {ord}

\DeclareMathOperator {\im} {Im}

\DeclareMathOperator {\SL} {SL}

\DeclareMathOperator {\GL} {GL}

\DeclareMathOperator {\seq} {\subseteq}
\DeclareMathOperator {\C} {\mathbb{C}}
\DeclareMathOperator {\R} {\mathbb{R}}
\DeclareMathOperator {\h} {\mathbb{H}}
\DeclareMathOperator {\Z} {\mathbb{Z}}
\DeclareMathOperator {\Q} {\mathbb{Q}}

\newcommand{\F}{\mathbb{F}}
\newcommand{\N}{\mathbb{N}}

\newcommand{\re}{\mathrm{Re}}

\DeclareMathOperator{\Jdeg}{deg_\J}
\DeclareMathOperator{\Jord}{ord_\J}

\newcommand{\J}{\mathbf{j}}
\newcommand{\Y}{\mathbf{Y}}

\theoremstyle {plain}
\newtheorem {theorem}{Theorem}[section]

\newtheorem {lemma} [theorem] {Lemma}
\newtheorem {proposition} [theorem] {Proposition}
\newtheorem {corollary} [theorem] {Corollary}

\newtheorem {claim} {Claim}[theorem]

\newtheorem{innercustomthm}{Theorem}
\newenvironment{customthm}[1]
  {\renewcommand\theinnercustomthm{#1}\innercustomthm}
  {\endinnercustomthm}

\theoremstyle {definition}

\newtheorem {definition}[theorem]{Definition}
\newtheorem{example} [theorem] {Example}
\newtheorem* {notation} {Notation}

\theoremstyle {remark}
\newtheorem {remark} [theorem] {Remark}
\numberwithin{equation}{section}

\newenvironment{claimproof}{\bgroup\begin{proof}}{\end{proof}\egroup}

\calclayout

\begin {document}

\title[Equations involving $j,j',j''$]{Equations involving the modular $j$-function and its derivatives}

\author{Vahagn Aslanyan}
\email{Vahagn.Aslanyan@manchester.ac.uk}
\address{Department of Mathematics, University of Manchester, Manchester, M13 9PL, UK}
\author{Sebastian Eterovi\'c}
\email{sebastian.eterovic@univie.ac.at}
\curraddr{Kurt G\"odel Research Center, Universit\"at Wien, 1090 Wien, Austria}

\author{Vincenzo Mantova}
\email{V.L.mantova@leeds.ac.uk}
\address{School of Mathematics, University of Leeds, Leeds, LS2 9JT, UK}

\date{}

\thanks{VA was Supported by Leverhulme Trust Early Career Fellowship ECF-2022-082 at the University of Leeds (where most of this work was done) and by EPSRC Fellowship EP/X009823/1 and DKO Fellowship at the University of Manchester. SE and VM were supported by EPSRC Fellowship EP/T018461/1 at the University of Leeds. For the purpose of open access, the authors have applied a Creative Commons Attribution (CC BY) licence to any Author Accepted Manuscript version arising from this submission.}

\keywords{Existential closedness, $j$-function, zero estimates, Rouché's theorem}

\subjclass[2020]{11F03, 11F23, 11U09}

\vspace*{-1cm}

\begin{abstract}
We show that for any polynomial $F(X,Y_0,Y_1,Y_2) \in \C[X, Y_0, Y_1, Y_2]$, the equation $F(z,j(z),j'(z),j''(z))=0$ has a Zariski dense set of solutions in the hypersurface $F(X,Y_0,Y_1,Y_2)=0$, unless $F$ is in $\C[X]$ or it is divisible by $Y_0$, $Y_0-1728$, or $Y_1$.

Our methods establish criteria for finding solutions to more general equations involving periodic functions. Furthermore, they produce a qualitative description of the distribution of these solutions.
\end{abstract}

\maketitle

\section{Introduction}\label{sec:introduction}

The problem of determining which (systems of) equations involving certain classical transcendental functions of a complex variable have solutions is a natural question at the intersection between complex geometry, model theory, and number theory. In complex geometry, it is a form of analytic Nullstellensatz for the given functions; in model theory, it plays an important role in the definability properties of the functions involved; and in number theory, it is related to Schanuel's conjecture and its analogues (given by special cases of the Grothendieck--Andr\'e generalised period conjecture). Often the function under consideration is of arithmetic importance. Examples of such classical functions are the exponential functions of semi-abelian varieties and Fuchsian automorphic functions. In this paper we focus on the modular $j$-function and its derivatives.

{The first conjecture in this area arose from Zilber's work on the model theory of complex exponentiation \cite{Zil2005,Zil2002,Zil2015}. It is now referred to as the \textit{Exponential (Algebraic) Closedness} conjecture or \textit{Zilber's Nullstellensatz}, and predicts when systems of equations involving addition, multiplication, and complex exponentiation have solutions in the complex numbers. We refer to the general version of the problem as \textit{Existential Closedness}, or EC for short. An EC conjecture for the $j$-function was proposed in \cite[\S 1]{AK2022}; in geometric terms it states that any algebraic variety $V\subseteq \C^{2n}$ satisfying geometric conditions known as \emph{freeness} and \emph{broadness}, intersects the $n$-fold graph of the $j$-function. The definition of these geometric conditions is long and will not be used in the present work, so we refer the interested reader to \cite[\S 2.2]{AK2022}, but, informally, freeness and broadness ensure that the equations defining $V$ do not break any functional properties of $j$ coming from the linear-fractional action of $\GL_2^+(\Q)$ (where $+$ denotes positive determinant) on the upper half-plane, as well as not contradicting a conjecture on transcendental values of the $j$-function analogous to Schanuel's conjecture for exponentiation (see \cite[Conj.~1.1]{AK2022} and \cite[\S 6.3]{AEK2023} for the statement of this conjecture).

If one somehow knows that an algebraic variety $V$ does intersect the graph of $j$, a very natural next question is to determine how these intersection points are distributed within $V$. For instance, one may ask whether these points are Zariski dense in $V$. We remark that if $V \subseteq \C^{2n}$ satisfies the above-mentioned geometric conditions of freeness and broadness, then for any Zariski open subset $V'\subseteq V$ it is possible to construct an algebraic variety $W \seq \C^{2(n+1)}$ which is also free and broad and projects onto $V'$. Thus, if we assume EC then, applying it to $W$, we deduce that $V'$ intersects the graph of $j$. Since $V'$ was an arbitrary Zariski open subset of $V$, we conclude that the intersection of $V$ with the graph of $j$ is Zariski dense in $V$. }

In the same work \cite{AK2022}, the authors also proposed an extension of the conjecture incorporating the derivatives of $j$ (see \cite[Conj. 1.6]{AK2022}). This version of EC is often referred to as \textit{Existential Closedness with Derivatives}, or ECD for short. This time the variety $V$ in question is a subset of $\mathbb{C}^{4n}$, and the conjecture states that if $V$ satisfies analogous geometric notions of freeness and broadness {(again related to a form of Schanuel's conjecture, now involving $j$ and its derivatives)}, then $V$ intersects the $n$-fold graph of the map $z\mapsto (j(z),j'(z),j''(z))$. For the definitions and precise statements of these conjectures, see \cite{AK2022,AEK2021,Asl2022}.
Note that we do not consider the third and higher derivatives of $j$ as these are rational over $j,j',j''$ (see (\ref{eq:j})).
As with EC, ECD implies that if $V\subseteq\C^{4n}$ satisfies the geometric conditions of freeness and broadness, then $V$ has a Zariski dense set of points of the desired form.

Very few cases of ECD have been proven, in comparison to EC where various families of varieties in $\C^{2n}$ have been shown to satisfy the conjecture.
Prior to the present work, only very special cases of ECD had been solved, proving solvability of some simple equations involving just $j'$ (so not combining it with $j$ or $j''$), see \cite{EH2021,Gal2021}. An ECD statement for ``blurrings'' (certain multi-valued twists) of $j$ was obtained in \cite{AK2022}. All of these papers mostly focused on, and established stronger results for, the EC conjecture for the $j$-function (without derivatives). These results are analogous to their exponential counterparts, namely, \cite{BM2017,DFT2018,AKM2023,Kir2019,Zil2002,Gal2023}. Although different methods have been used across these works, one common feature is that they all exploit in some way the periodicity of $\exp$ or the $\SL_2(\Z)$-invariance of $j$.
Incorporating the derivatives of $j$ into the equations presents then a significant new challenge, as $j'$ and $j''$ are no longer $\SL_2(\Z)$-invariant.
It is also worth noting that a differential version of ECD was obtained in \cite{AEK2021}, and that it is so far the only setting where a full Existential Closedness statement is proved for $j$ with derivatives, and the same method also applies to $\exp$.

In this article we prove the ECD conjecture when $n=1$, which precisely states that any algebraic variety $V\subseteq\mathbb{C}^4$ of dimension 3 without constant coordinates contains (a Zariski dense set of) points of the form $(z,j(z),j'(z),j''(z))$.
This amounts to checking exactly which equations of one complex variable involving only $z,j(z),j'(z), j''(z)$ have solutions, and whether these solutions are Zariski dense. Note that for $n=1$ broadness of $V$ just means $\dim V \geq 3$ hence the only non-trivial case is $\dim V = 3$. On the other hand, freeness for $n=1$ means $V$ has no constant coordinates. Nevertheless, we will even be able to decide what happens when $V$ does have a constant coordinate.

Our first main result establishes the existence of solutions in all non-trivial cases.

\begin{theorem}\label{thm: intro existence of sol - not dense}
Let $F(X,Y_0,Y_1,Y_2) \in \C[X,Y_0,Y_1,Y_2] \setminus \C[X]$. Then the equation $F(z,j(z),j'(z),j''(z))=0$ has infinitely many solutions.
\end{theorem}

The proof of \Cref{thm: intro existence of sol - not dense} is based on a generalisation of the methods of \cite{EH2021}, which use Rouch\'e's theorem from classical complex analysis to establish some cases of EC for $j$ (without derivatives). \Cref{thm: intro existence of sol - not dense} can be seen as an analogue of the classical fact that every irreducible polynomial $p(X,Y)\in\C[X,Y]$ which depends on $Y$ has infinitely many zeroes of the form $(z,\exp(z))$, unless $p = cY$ for some $c\in\C^\times$.

Throughout the paper, all algebraic subvarieties of $\C^4$ will be defined by polynomials in the ring $\C[X, Y_0,Y_1,Y_2]$.

Our main goal is to obtain a much stronger version of \Cref{thm: intro existence of sol - not dense}. We show that for any polynomial $F(X,Y_0,Y_1,Y_2)$ the set \[ \{ (z,j(z),j'(z),j''(z)) \in \h\times \C^3: F(z,j(z),j'(z),j''(z))=0 \} \] is Zariski dense in the hypersurface $F(X,Y_0,Y_1,Y_2)=0$, unless $F$ is divisible by an explicit (finite) list of polynomials. In this case we say that the equation $F(z,j(z),j'(z),j''(z))=0$ has a \textit{Zariski dense set of solutions} (see \Cref{def:zariski-dense}), that is, by a solution of such an equation we understand a tuple $(z_0,j(z_0),j'(z_0),j''(z_0))$ rather than just $z_0$.

{We remind the reader that this is equivalent to establishing certain cases of ECD for subvarieties of $\C^8$: given a hypersurface $V \subseteq \C^4$ and a Zariski open dense $V' \subseteq V$, there is $W \subseteq \C^8$ free and broad which projects onto $V'$, such that $W$ intersects the graph of $j$ and its derivatives if and only if $V'$ does. For instance, the system $\{ j''(z)=0,~j(z)\neq 0 \}$ has a solution if and only if $\{ j''(z_1)=0,~j(z_1)z_2=1 \}$ does.}

The bulk of the paper is focused on proving the Zariski density of the set of solutions, which the proof of \Cref{thm: intro existence of sol - not dense} does not provide. For instance, the solutions of the equation $zj''(z)+(z^3+1)j'(z)^2+j'(z)j(z)^7=0$ found via the proof Theorem~\ref{thm: intro existence of sol - not dense} are the $\SL_2(\Z)$-conjugates of $\rho = -\frac{1}{2}+\frac{\sqrt{3}}{2}i$. These are obviously not Zariski dense, for it is well known that {$j(\gamma\rho)=j'(\gamma\rho)=j''(\gamma\rho)=0$ for every $\gamma \in \SL_2(\Z)$} (see \cite[p.~40]{lang-ellipticfunctions}). Indeed, {Zariski density requires at least that the solutions are not contained in finitely many $\SL_2(\Z)$-orbits, except for when the equation is of the form $\prod_k (j(z) - u_k) = 0$ for some $u_k \in \C$}.

{The zeroes of $j'$ are in fact problematic:} observe that
\[ \forall z \in \h \left[ j(z)(j(z)-1728)=0 \Longleftrightarrow j'(z)=0 \right]. \]
This immediately gives that the three equations $j(z)=0$, $j(z)-1728=0$, and $j'(z)=0$ do not have Zariski dense sets of solutions. Our main result shows that these are essentially the only non-examples.

\begin{theorem}\label{thm: intro - main thm}
        For any polynomial $F(X,Y_0,Y_1,Y_2) \in \C[X, Y_0, Y_1, Y_2]\setminus \C[X]$ which is coprime to $Y_0(Y_0-1728)Y_1$, the equation $F(z,j(z),j'(z),j''(z))=0$ has a Zariski dense set of solutions, i.e.\ the set \[ \{ (z,j(z),j'(z),j''(z)) \in \h\times \C^3: F(z,j(z),j'(z),j''(z))=0 \}\] is Zariski dense in the hypersurface $F(X,Y_0,Y_1,Y_2)=0$.
\end{theorem}

\begin{remark}
    \Cref{thm: intro - main thm} implies that for every rational function $G(X,Y_0,Y_1,Y_2) \in \C(X,Y_0,Y_1,Y_2)$, the function $G(z,j(z),j'(z),j''(z))$ has a zero unless $G$ is of the form
    \[ \frac{Y_0^s(Y_0-1728)^tY_1^\ell}{H(X,Y_0,Y_1,Y_2)} \]
    where $H$ is a polynomial and $s, t, \ell \in \mathbb{N}$.
\end{remark}

A special case of \Cref{thm: intro - main thm} is that the equation $j''(z)=0$ has a Zariski dense set of solutions.\footnote{In particular, the ramification points of $j'$ are not contained in finitely many $\SL_2(\Z)$-orbits.} In this case, even proving that there is a solution outside the $\SL_2(\Z)$-orbit of $\rho$ is highly non-trivial, see \S\ref{subsec:homeq(j,j'j,'')}.

\begin{remark}
   In \cite{Ete2022}, the author studies the problem of finding \emph{generic} solutions to equations involving $j$ (and its derivatives) under the assumption that the system has a Zariski dense set of solutions.
    In particular, combining \Cref{thm: intro - main thm} with \cite[Theorem 6.5]{Ete2022} we get that there is a countable field $C_j\subseteq \C$ such that for any irreducible hypersurface $V\subset \C^4$ satisfying the conditions of \Cref{thm: intro - main thm}, if $V$ is not definable over $C_j$, then for any finitely generated subfield $K\subset \C$ over which $V$ can be defined there is a point of the form $(z,j(z),j'(z),j''(z))\in V$ such that
    \begin{equation*}
        \mathrm{tr.deg.}_K K(z,j(z),j'(z).j''(z)) = \dim V = 3.
    \end{equation*}
\end{remark}

To prove Theorem~\ref{thm: intro - main thm} we establish general criteria for the solvability of certain equations involving periodic functions (see \S\ref{sec:periodic}). The following proposition is a special case of those criteria.

\begin{definition}\label{def: merom at infty}
    {A meromorphic function $f:\h\to\C$ is \emph{1-periodic} if $f(z + 1) = f(z)$ for every $z \in \h$. 
    Every such function induces a meromorphic function $\tilde{f}(q)$ on the punctured unit disc by performing the change of variable $q=\exp(2\pi iz)$}. We say that $f$ is \emph{meromorphic at $i\infty$} if $\tilde{f}$ is meromorphic at $0$.
\end{definition}

\begin{proposition}\label{prop: intro periodic pole}
Let $f_0,\ldots, f_{n} : \h \to \C$ be $1$-periodic functions {which are meromorphic on $\h \cup \{i\infty\}$}.    Suppose that for some $k$ one of the following conditions is satisfied:
\begin{itemize}
    \item there is $\tau\in\h$ such that $\frac{f_k}{f_n}(z) \to \infty$ as $z \to \tau \in \h$, or
    \item $\frac{f_k}{f_n}(z) \to \infty$ as $\im(z)\to+\infty$.
\end{itemize}
Then there is a sequence of points $\{z_m\}_{m\in\N}\subseteq \h$ with $z_m\neq \tau$ and $z_m\to \tau$ in the first case, or $\im(z_m)\to+\infty$ and $0\leq \re(z_m)\leq 2$ in the second case, such that for all sufficiently large $m$ the point $z_m+m$ is a solution to the equation
\[f_n(z)z^n+f_{n-1}(z)z^{n-1}+\ldots+f_0(z) = 0.\]
\end{proposition}

Let us consider an example illustrating how we apply \Cref{prop: intro periodic pole} in practice. It also gives an idea of our approach in the general case.

\begin{example}
    Consider the equation
    \begin{equation}\label{eqn: intro example pole}
        j'(z)^2 + p(j(z)) = 0
    \end{equation}
    where either $p(j(z)) = j(z)(j(z)-1728)$ or $p(j(z)) = j(z)^2(j(z)-1728)$. First, we want to get an equivalent equation which is written as a sum of powers of $z$ with periodic coefficients. To that end we apply the $\SL_2(\Z)$-transformation $z \mapsto -\frac{1}{z}$ and, using the identities $j\left( -\frac{1}{z} \right) = j(z),~ j'\left( -\frac{1}{z} \right) = z^2j'(z)$, we get
    \begin{equation}\label{eqn: intro example pole 2}
        z^4 j'(z)^2 + p(j(z)) = 0.
    \end{equation}
    Thus we obtain an equation in a suitable form for using Proposition \ref{prop: intro periodic pole}, where $f_4 = (j')^2$, $f_3=f_2=f_1=0$ and $f_0 = p(j)$.
    When $p(j) = j(j-1728)$, the ratio $\frac{f_0}{f_4}=\frac{j(j-1728)}{(j')^2}$ has a pole at $\tau = \rho$. When $p(j) = j^2(j-1728)$, the ratio $\frac{f_0}{f_4}=\frac{j^2(j-1728)}{(j')^2}$ has no finite poles, but it has limit $\infty$ as $\im(z) \to  +\infty$.

    Thus, by \Cref{prop: intro periodic pole}, there is a sequence $z_m$ with $z_m\to \tau$ and $z_m\neq \tau$ in the first case, or $\im(z_m)\to +\infty$ and $0\leq \re(z_m)\leq 2$ in the second case, such that  for all sufficiently large $m$ the point $z_m+m$ is a solution to the equation \eqref{eqn: intro example pole 2}. This already implies that the solutions of  \eqref{eqn: intro example pole 2} intersect infinitely many $\SL_2(\Z)$-orbits.

    To deduce Zariski density of these solutions, suppose that all of them are also solutions of another independent equation $G(z,j(z),j'(z),j''(z))=0$. Combining this and \eqref{eqn: intro example pole 2} we can eliminate $z$ and end up with an equation $H(j,j',j'')=0$. Now, our assumption means that $H(j,j',j'')$ vanishes at $z_m+m$, hence also at $z_m$ by periodicity. This is not possible, for a 1-periodic holomorphic function, meromorphic at $i\infty$, cannot have infinitely many zeroes with real part bounded from above and below and imaginary part bounded from below.
    This then implies the Zariski density of solutions of \eqref{eqn: intro example pole}.
\end{example}

We also note that our criteria can be applied to more general periodic functions, beyond polynomials of $j,j',j''$, such as $\exp$ or the Weierstrass $\wp$-function. For instance, \Cref{prop: intro periodic pole} implies that the function $j'(z)z+\exp(2\pi i z)$ has infinitely many zeroes around the points $i+m$ where $m$ is a large integer.

\subsection{Structure of the paper}
\begin{description}
    \item[\S \ref{sec:prelim}] We go over some basic preliminaries about the $j$-function and its derivatives. We also give the definition of Zariski density used in Theorem \ref{thm: intro - main thm}.
    \item[\S \ref{sec:existence}] Here we prove Theorem \ref{thm: intro existence of sol - not dense} by extending the methods of \cite{EH2021}, which are based on Rouch\'e's theorem.
    \item[\S\ref{sec:periodic}] We prove criteria for the existence and distribution of solutions of equations involving periodic functions, which combined imply Proposition~\ref{prop: intro periodic pole} (but are significantly more general). The approach used here involves Rouch\'e's theorem, the Argument Principle, and some elementary methods from valuation theory. These methods do not appear in later sections of the paper.
    \item[\S \ref{sec:zariskidensity}]
    We use the results of the previous section to obtain concrete criteria for proving Zariski density of equations of the form $F(z,j(z),j'(z),j''(z))=0$. These criteria are about the presence of poles in quotients of certain polynomials only in $j,j',j''$.
    \item[\S\ref{sec:zeroesetimates}] We produce zero estimates for polynomials in $z$, $j$, $j'$, $j''$ in order to determine when the quotients mentioned above have poles.
    \item[\S\ref{sec:eq(z,j,j',j'')}] First, in \S\ref{subsec:homeq(j,j'j,'')} we prove Zariski density for \emph{$\J$-homogeneous} equations (\Cref{def:j-homogeneous}), which only involve $j(z)$, $j'(z)$ and $j''(z)$, including the equation $j''(z)=0$. Finally, in \S\ref{subsec:proofmainthm} we prove Theorem \ref{thm: intro - main thm} in full generality.
\end{description}

\section{Preliminaries}
\label{sec:prelim}
Let $\h$ denote the complex upper half-plane $\{z\in \mathbb{C}:\mathrm{Im}(z)>0\}$. The group $\mathrm{GL}_2^+(\R)$ of $2\times 2$ real invertible matrices with positive determinant acts on $\h$ via linear fractional transformations:
\[gz:=\frac{az+b}{cz+d} \quad \text{for} \quad g=\begin{pmatrix}a & b\\ c & d\end{pmatrix}\in\mathrm{GL}_2^+(\R).\]
This action can be seen as a restriction of the action of $\mathrm{GL}_2(\C)$ on the Riemann sphere $\C\cup\{\infty\}$. The \emph{modular group} is defined as
\[\mathrm{SL}_2(\Z):=\left\{\begin{pmatrix}a & b\\ c & d\end{pmatrix}\in \mathrm{GL}_2^+(\R): a,b,c,d\in\Z \text{ and } ad-bc=1\right\}.\]
As a group, $\mathrm{SL}_2(\Z)$ is generated by two elements: $\begin{pmatrix}
1 & 1 \\
0 & 1
\end{pmatrix}$ and $\begin{pmatrix}
0 & -1 \\
1 & 0
\end{pmatrix}$, which correspond to the actions $z\mapsto z+1$ and $z\mapsto -\frac{1}{z}$, respectively.

The \emph{modular $j$-function} is defined as the unique $\mathrm{SL}_2(\Z)$-automorphic function $j:\h\to \C$ satisfying $j(\rho)=0$ (recall that $\rho:=\exp\left(\frac{2\pi i}{3}\right)$; this notation will be kept throughout the paper), $j(i)=1728$ and $j(\infty)=\infty$ (this last condition should be understood as $\lim_{\mathrm{Im}(z)\to+\infty}j(z)=\infty$). In particular, this means that $j$ satisfies
\[j(\gamma z)=j(z) \text{ for every }\gamma \text{ in } \mathrm{SL}_2(\Z) \text{ and every }z\text{ in }\h,\]
and it is in particular 1-periodic, and by assumption meromorphic at $i\infty$. Its Fourier expansion (also known as a $q$-expansion) is of the form
\begin{equation}\label{eq:j-fourier-expansion}
j(z)=q^{-1}+744+\sum_{k=1}^{\infty}a_kq^k,\, \text{ with }q:=\exp(2\pi i z)  \text{ and }a_k\in \C.
\end{equation}
In fact, $a_k\in\Z$ for every $k\in\N$. The $j$-function induces an analytic isomorphism of Riemann surfaces $\mathrm{SL}_2(\Z) \backslash \h\simeq \C$ (see \cite[Chapter 3, \S3]{lang-ellipticfunctions}).

\begin{figure}
    \centering
    \includegraphics{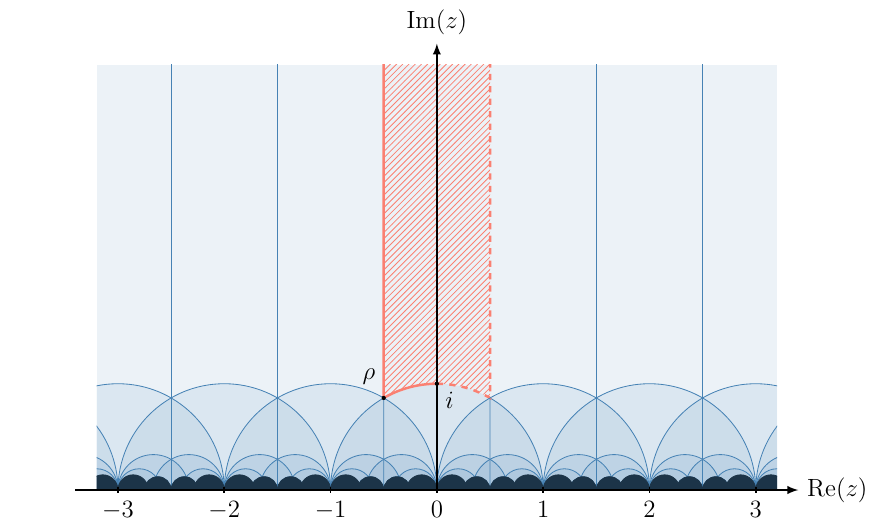}
    \caption{\label{fundamental-domain}The fundamental domains of the action by $\SL_2(\Z)$, where $\F$ is highlighted by the striped background.}
\end{figure}

Since $j$ is invariant under the action of $\mathrm{SL}_2(\Z)$, we can study the behaviour of $j$ by looking at fundamental domains of the action of $\mathrm{SL}_2(\Z)$ on $\h$. The \emph{standard fundamental domain} is the set
\begin{equation*}
    \F:=\left\{z\in\C: -\frac{1}{2}\leq\mathrm{Re}(z)<\frac{1}{2},\, |z|\geq 1,\,\left( |z|=1\implies -\frac{1}{2}\leq \mathrm{Re}(z)\leq 0 \right) \right\}.
\end{equation*}
We let $\overline{\F}$ denote the Euclidean closure of $\F$ (within the Riemann sphere). A diagram of the standard fundamental domain along with some of its $\mathrm{SL}_2(\Z)$-translates is given in \Cref{fundamental-domain}.
When we refer to a \emph{fundamental domain}, we will always mean a set of the form $\gamma\F$ for some $\gamma\in\SL_2(\Z)$.

It is well known that $j$ satisfies the following third-order differential equation (and none of lower order, \cite{mahler}):\footnote{Observe that the denominators in the equation correspond to the polynomials that must be omitted in \Cref{thm: intro - main thm}.}
\begin{equation}
    \label{eq:j}
    0 = \frac{j'''}{j'} - \frac{3}{2}\left(\frac{j''}{j'}\right)^{2} + \frac{j^{2} -1968j + 2654208}{2j^{2}(j-1728)^{2}}\left(j'\right)^{2}.
\end{equation}
This shows that the derivatives of $j$ of order at least 3 are rational over $j$, $j'$ and $j''$. Mahler's result \cite{mahler} implies that $j$, $j'$ and $j''$ are algebraically independent over $\C$.

The functions $j$, $j'$ and $j''$ are all 1-periodic and meromorphic at $i\infty$, and by differentiating (\ref{eq:j-fourier-expansion}) we can obtain the $q$-expansions of $j'$ and $j''$:
\begin{align*}
    j'(z) &= -\frac{2\pi i}{q} + 2\pi i\sum_{k\geq 1}a_k kq^{k}\\
    j''(z)&= -\frac{4\pi^2}{q} -4\pi^2 \sum_{k\geq 1}a_k k^2q^{k}.
\end{align*}

Observe that $\Q\cup\{\infty\}$ forms a single orbit under the action of $\mathrm{SL}_2(\Z)$. We call these elements the \emph{cusps} of $j$. Given a fundamental domain $\gamma\F$, the \emph{cusp of $\gamma\F$} is the unique element of $\Q\cup\{\infty\}$ contained in the Euclidean closure of $\gamma\F$ (where the closure is taken in the Riemann sphere).

Using (\ref{eq:j-fourier-expansion}) and the $q$-expansions of $j'$ and $j''$, it is easy to see that for every $x\in \R$ we have that each of the expressions $j(x+iy)$, $j'(x+iy)$ and $j''(x+iy)$ grows exponentially to $\infty$ as $y\to+\infty$. In this case we sometimes write $z\to i\infty$ to emphasise that $z$ approaches $\infty$ by increasing its imaginary part, while the real part remains bounded. A similar behaviour takes place when $z$ approaches a rational number from within a fixed fundamental domain containing that rational number in its euclidean closure. We will give a precise description of this behaviour in \S\ref{sec:periodic}.

We finish this section with the definition of what we mean by finding a Zariski dense set of solutions to an equation.

\begin{definition}\label{def:zariski-dense}
    Let $F(X,Y_0,Y_1,Y_2)$ be a polynomial over $\C$. We say the equation \[ F(z,j(z),j'(z),j''(z))=0 \] has a \textit{Zariski dense set of solutions} if for any polynomial $G(X,Y_0,Y_1,Y_2)$ which is not divisible by some
    irreducible factor of $F$, there is $z_0\in \h$ such that \[ F(z_0,j(z_0),j'(z_0),j''(z_0))=0 \text{ and }G(z_0,j(z_0),j'(z_0),j''(z_0))\neq 0.\]
\end{definition}

Clearly, it suffices to prove Zariski density for irreducible polynomials to obtain Theorem~\ref{thm: intro - main thm}, so from now on we will reduce to the case where $F$ is irreducible.

\section{Existence of Solutions}
\label{sec:existence}
We start with the \emph{Rouch\'e method} for proving the existence of solutions, but not yet their Zariski density.
We first recall the crucial theorem.

\begin{theorem}[Rouch\'e, see e.g.~{{\cite[Chapter VI, \S1, Theorem 1.6]{lang:complexanalysis}}}]
    \label{thm:rouche} Let $f,g$ be meromorphic functions on a complex domain $\Omega$. Let $C$ denote a simple closed curve which is homologous to $0$ in $\Omega$ and such that $f$ has no zeroes or poles on $C$. If the inequality
    \[|g(z)|< |f(z)|\]
    holds for all $z$ on $C$, then the difference between the numbers of zeroes and poles in the interior of $C$ for the functions $f+g$ and $f$ is the same.
\end{theorem}

{It is well known that the Euclidean closure of the $\SL_2(\Z)$-orbit of any point in $\h$ (where the closure is taken within the Riemann sphere) only accumulates at the boundary of $\h$, that is, at $\R\cup\{\infty\}$.
The following lemma will help us choose convenient sequences within any given orbit converging to points in $\R$}.

\begin{lemma}\label{lemma-gamma-bounded}
Let $z \in \h$ and $u \in \R$.
\begin{enumerate}[label=(\roman*)]
    \item If $ u \in \R \setminus \Q$ and the sequence $\gamma_k=\begin{pmatrix}
a_k & b_k \\
c_k & d_k
\end{pmatrix} \in \SL_2(\mathbb{Z})$ is such that $\gamma_k z \to u$ as $k \to +\infty$, then { $|c_k|\to +\infty$ and $\frac{a_k}{c_k}\to u$ as $k \to +\infty$.}

\item\label{rational-cusp} If $u = \frac{a}{c} \in \Q$ with $\gcd(a,c)=1$, then there is a sequence $\gamma_k=\begin{pmatrix}
a & b_k \\
c & d_k
\end{pmatrix} \in \SL_2(\mathbb{Z})$ such that $|b_k|, |d_k| \to +\infty$ and $\gamma_k z \to u$ as $k \to +\infty$.
\end{enumerate}
\end{lemma}
\begin{proof}
\noindent
\begin{enumerate}[label=(\roman*)]
 \item  {We first show that $|c_k|\to+\infty $.}
 Since every subsequence of $\gamma_k z$ tends to $u$, it suffices to show that $c_k$ is unbounded. Assume it is bounded, then we can choose a subsequence where the value of $c_k$ is constant. So assume $c_k=c$ is a constant sequence.

Conversely, assume now that $\gamma_kz\to u$.
Let $z = x+iy$. If $d_k$ is also bounded then we may assume it is constant. Then $a_kz+b_k = (a_kx+b_k) + a_k yi$ must be convergent, hence $a_k$ must be convergent, and so constant. Then $b_k$ is also constant, for $a_kd_k-b_kc_k=1$, a contradiction.

Thus, we may assume $|d_k| \to +\infty$. Then $$\frac{\frac{a_k}{d_k}z+\frac{b_k}{d_k}}{\frac{c}{d_k}z+1}\to u.$$ Therefore, $$\frac{a_k}{d_k}z+\frac{b_k}{d_k} = \left( \frac{a_k}{d_k} x + \frac{b_k}{d_k}\right) + \frac{a_k}{d_k} y i \rightarrow u.$$
This implies $$\frac{a_k}{d_k}\to 0,~ \frac{b_k}{d_k} \to u.$$
On the other hand, $a_kd_k-b_kc=1$, hence $a_k - \frac{b_k}{d_k}c = \frac{1}{d_k}\to 0$. Thus, $a_k \to uc$ which means $a_k = a \in \Z$ is constant. But then $u = \frac{a}{c} \in \Q$.

{Since $|c_k|\to+\infty$, then using that
\[\frac{a_kz+b_k}{c_kz+d_k}\cdot\frac{c_k}{a_k} = \frac{a_kz + b_k}{a_kz + b_k + \frac{1}{c_k}}\to 1\]
we see that $\gamma_k z$ and $\frac{a_k}{c_k}$ must have the same limit.}

\item Since $u = \frac{a}{c}$ with $\gcd(a,c)=1$, there are integers $m, l $ such that $am-cl=1$. Choose $b_k = l + ka,~ d_k = m + kc$. Then
\[ \lim_{k \to +\infty}\frac{az+b_k}{cz+d_k} = \lim_{k \to +\infty}\frac{\frac{a}{d_k}z+\frac{b_k}{d_k}}{\frac{c}{d_k}z+1}= \lim_{k \to +\infty} \frac{b_k}{d_k} =  \lim_{k \to +\infty} \frac{l + ak}{m+ck} = \frac{a}{c} = u.\qedhere \]
\end{enumerate}
\end{proof}

In order to ease notation, we will start using bold-faced letters to denote vectors, so we set $\Y:=(Y_0,Y_1,Y_2)$ and $\J:=(j,j',j''): \h \to \C^3$.

We are now ready to prove Theorem \ref{thm: intro existence of sol - not dense} which we restate below for convenience.

\begin{customthm}{\ref{thm: intro existence of sol - not dense}}
For every $F(X,\Y) \in \C[X,\Y] \setminus \C[X]$ the equation $F(z,\J(z))=0$ has infinitely many solutions.
\end{customthm}

\begin{proof}
If $F$ does not depend on $Y_1$ and $Y_2$ then we are done by the results of \cite{EH2021}. So assume $F$ depends on $Y_1$ or $Y_2$. The argument below is a generalisation of the method of \cite{EH2021}.

Let $r(X)\coloneqq -F(X,0,0,0)$ and $G(X,\Y)\coloneqq  F(X,\Y)+r(X)$. Further let $f(z)\coloneqq G(z,\J(z))$. Then we want to solve the equation $$f(z) = r(z).$$
Notice that $f(\rho) = 0$, for $j(\rho) = j'(\rho) = j''(\rho)=0$. Let $B\subseteq \C$ be a closed disc centred at $\rho$ with sufficiently small radius such that
 $j'(z)\neq 0$ and $j''(z)\neq 0$ for $z\in  B\setminus \{ \rho \}$.

Pick a point $u = \frac{a}{c} \in \Q$, and choose a sequence $\gamma_k \in \SL_2(\Z)$ as in \Cref{lemma-gamma-bounded}\eqref{rational-cusp} such that $\gamma_k z \to u$ as $k \to +\infty$ for any $z \in \h$. Let $B_k\coloneqq  \gamma_k B$. By compactness of $B$, the function $r(\gamma_kz)$ tends to $r(u)$ uniformly for $z\in B$.

We have
\begin{align*}
    f(\gamma_k z) &=  G\left(\gamma_kz, j(z), (cz+d_k)^2 j'(z), (cz+d_k)^4 j''(z) + 2c(cz+d_k)^3 j'(z)\right)\\
    &= G\left(\gamma_kz, j(z), d_k^2\left(\tfrac{c}{d_k}z+1\right)^2 j'(z), d_k^4\left(\left(\tfrac{c}{d_k}z+1 \right)^4 j''(z) + 2\tfrac{c}{d_k}\left(\tfrac{c}{d_k}z+1\right)^3 j'(z)\right)\right).
\end{align*}

Consider the polynomial $G(X,Y_0,T^2Y_1,T^4Y_2)$ as a polynomial of $T$. It clearly has positive degree, for otherwise $G$ (and hence $F$) would not depend on $Y_1$ nor $Y_2$. Let its leading term be $H(X,\Y)\cdot T^m$. Since $\frac{c}{d_k}\to 0$, we see that
\[ f(\gamma_k z) = d_k^m \cdot H(u,\J(z)) + o(d_k^m) \text{ as } k \to +\infty. \]

We can now shrink $B$ to make sure that $H(u,\J(z))\neq 0$ on $\partial B$, so it is uniformly bounded away from $0$ for $z \in \partial B$. In particular, $f(\gamma_k z)$ approaches infinity as $k \to +\infty$ uniformly for $z \in \partial B$. So for sufficiently large $k$ the inequality $|f(z)| > |r(z)|$ holds for all $z \in \partial B_k = \partial(\gamma_k B) = \gamma_k \partial B$, and we can apply Rouch\'e's theorem to these functions. Since $f$ has a zero in $B_k$, namely $\gamma_k \rho$, so does $f-r$.
\end{proof}

\begin{remark}
The following more general statement can be proven by the same argument.

\emph{Let $F(X,\Y) \in \C[X,\Y] \setminus \C[X]$. Let $U \subseteq \C$ be an open set such that $U \cap \R \neq \emptyset$ and let $f:U \rightarrow \C$ be a holomorphic function.  Then the equation $F(z,\J(z))=f(z)$ has infinitely many solutions.}
\end{remark}

As mentioned in \S\ref{sec:introduction}, the proof of Theorem~\ref{thm: intro existence of sol - not dense} does not guarantee a Zariski dense set of solutions. For instance, if $F(X,0,0,0)\equiv 0$, that is, $F$ has no term depending only on $X$, then the only solutions yielded by the method above are the $\SL_2(\Z)$-conjugates of $\rho$. In order to establish Zariski density we will look at a refinement of the procedure, where we transform the equation by convenient elements of $\SL_2(\Z)$. This will be done starting in \S\ref{sec:zeroesetimates}, but first, in the next section, we will develop some tools to study equations involving periodic functions.

\section{Solvability of Certain Equations Involving Periodic Functions}
\label{sec:periodic}

In this section we establish some general criteria for the solvability of equations involving periodic functions and, in particular, prove Proposition \ref{prop: intro periodic pole}. We remark that this section is independent in many ways from the rest of the paper as the results we prove make no reference to $j$ or its derivatives, and in particular the methods developed here will not reappear in the following sections.

We recall that given a meromorphic function $f$, not identically 0, and a point $z_0$, the \emph{order} of $f$ at $z_0$ is the unique integer $n$ such that $(z-z_0)^{-n}f(z)$ is holomorphic and non-zero at $z_0$.

\begin{proposition}\label{prop: periodic finite pole}
Let $f_0,\ldots, f_{n} : \h \to \C$ be $1$-periodic meromorphic functions and let $-\ell$ be the minimum order of $\frac{f_k}{f_n}$ at a fixed $z_0 \in \h$ for $k = 0, \dots, n-1$.

If $\ell > 0$, then for any sufficiently small disc $D$ centred at $z_0$ and for every sufficiently large $m \in \Z$,
the  equation
\[f_n(z)z^n+f_{n-1}(z)z^{n-1}+\dots+f_0(z) = 0\]
has $\ell$ solutions, counted with multiplicity, in $m + (D \setminus \{z_0\})$.
\end{proposition}
\begin{proof}
    For simplicity, assume that $\frac{f_0}{f_n}$ has a pole at $z_0\in \h$ of order $\ell > 0$; the same proof will work in the general case with trivial modifications.

    Under the above assumptions, $\frac{f_n}{f_0}(z_0) = 0$, and moreover $\frac{f_k}{f_0}(z_0) \neq \infty$ for all $k$. Let $F(z)\coloneqq f_n(z)z^n+f_{n-1}(z)z^{n-1}+\dots+f_0(z)$. Consider the functions
    \[ G(z)\coloneqq \frac{f_n(z)}{f_0(z)}z^n \mbox{ and } H(z)\coloneqq \frac{f_{n-1}(z)}{f_0(z)}z^{n-1}+\dots+\frac{f_{1}(z)}{f_0(z)}z +1.\]
    Pick a small closed disc $D$ centred at $z_0$ such that the $f_k$'s have no zeroes nor poles in $D \setminus \{ z_0\}$. Since $\frac{f_k(z)}{f_0(z)}$ are periodic and bounded on $D$, for large enough $m$ we have
    \[ |G(z+m)| = \left| \frac{f_n(z)}{f_0(z)} \right| |z+m|^n >|H(z+m)| \text{ for } z\in \partial D.\]

    By Rouch\'e's \Cref{thm:rouche}, the number of zeroes of the functions $G(z+m)$ and $G(z+m)+H(z+m) ={f_0(z)}^{-1} F(z+m)$ inside $D$ is the same. Since the former has a zero at $z_0$ of order $\ell$ and no other zero, the latter must also have $\ell$ zeroes in $D$, counted with multiplicity. Thus, ${f_0(z)}^{-1}F(z+m)$ has $\ell$ zeroes in $D$, and so $f_0^{-1}F$ has $\ell$ zeroes in $m + D$, counted with multiplicity.

    Finally, note that $G(z_0 + m) + H(z_0 + m) = 0$ holds for at most $n-1$ values of $m$, thus for $m$ sufficiently large, the above $\ell$ solutions in $m + D$ are actually in $m + (D \setminus \{z_0\})$.

    This finishes the proof of the proposition when $\frac{f_0}{f_n}$ has a pole at $z_0$ of order greater than or equal to that of $\frac{f_k}{f_n}$ for $k = 1, \dots, n-1$. For when the maximum order of pole at $z_0$ is attained by $\frac{f_k}{f_n}$ for some $k \neq 0$, simply divide by $f_k(z)$ rather than $f_0(z)$ when defining $G$ and $H$.
\end{proof}

Following the notation of the proposition, when the $f_k$'s are polynomials in $j$, $j'$, $j''$ and some $\frac{f_k}{f_n}$ has a pole in $\h$, the above proposition applies. Instead, when the functions $\frac{f_k}{f_n}$ have no poles in $\h$, we will rely on the asymptotic behaviour of $\frac{f_k}{f_n}$ towards the boundary of $\h$. Specifically, we will prove an analogue of Proposition~\ref{prop: periodic finite pole} in the case where $\frac{f_k(z)}{f_n(z)} \to \infty$ as $z \to i\infty$.

As the following example shows, we may also need to consider non-periodic functions which are \emph{asymptotically periodic}. Dealing with those functions requires a considerably more sophisticated setup than the one in \Cref{prop: intro periodic pole}, so we first discuss the example in detail to clarify the choices made in the rest of this section.

\begin{example}

    Consider the equation \[ (j')^5 + j^2(j-1728)^2(j'')^2 + \alpha j^2(j-1728)(j')^3 = 0 \] where $\alpha$ is to be determined later.
    In order to write this equation as a polynomial in $z$ with periodic coefficients, we apply the $z\mapsto -\frac{1}{z}$ transformation\footnote{    {The action of the group $\SL_2(\Z)$ is generated by the transformations $z\mapsto z+1$ and $z\mapsto -\frac{1}{z}$. Since our functions are invariant under the former, it is natural to apply the latter; while $j$ is invariant under it, $j'$ and $j''$ are not, and we take advantage of this fact.}} and get
    \begin{align*}
    z^{10} (j'(z))^5 &+ z^8 j(z)^2(j(z)-1728)^2(j''(z))^2 + z^7 4j(z)^2(j(z)-1728)^2j'(z)j''(z)\\
    &+ z^6 (4j(z)^2(j(z)-1728)^2(j'(z))^2 + \alpha j(z)^2(j(z)-1728)(j'(z))^3) = 0.
    \end{align*}

    In this example the ratio of the coefficients of $z^8$ and $z^{10}$ actually has a pole at $i$. However, checking for poles among such ratios in a general equation requires sufficiently precise zero estimates at the conjugates of $\rho$ and $i$, which are hard to produce for polynomials involving $j''$
    (see \Cref{example: no zero estimate for j''}); on the other hand, we can provide sharp zero estimates for polynomials in $j$, $j'$ only (see \S\ref{sec:zeroesetimates}). The latter estimates turn out to be enough: for instance, when the original equation does not contain $z$, after the $z\mapsto -\frac{1}{z}$ transformation the coefficient of the lowest power of $z$ does not depend on $j''$. This is exemplified here by the coefficient of $z^6$.
    Hence our strategy hinges on the fact that the ratio between a particular coefficient not involving $j''$, which we can procure in all cases, and the leading coefficient has a pole or exponential growth in some fundamental domain.

    In this particular example, the ratio in question is between the coefficients of $z^6$ and $z^{10}$, thus the function
    \[ f(z):=\frac{4j^2(j-1728)^2 + \alpha j^2(j-1728)j'}{(j')^3}.\]
    We claim that $f(z)$ has no pole in $\h$. Indeed, easy calculations show that the orders of the numerator at $\rho$ and $i$ (and their $\SL_2(\Z)$-orbits) are equal to 6 and 3 respectively. The denominator has the same orders at these points, so $f$ has no poles. Moreover, choosing $\alpha = \frac{2}{\pi i}$ ensures the leading terms in the $q$-expansions of the two terms in the numerator cancel out. This then means that $f(z)$ tends to a constant as $z \to i\infty$. Therefore, \Cref{prop: intro periodic pole} cannot be applied in this situation.
    However, $f(z)$ has exponential growth as we approach $0$ from within a fundamental domain with a cusp at $0$ (in fact, we shall prove that $f(z)$ must have exponential growth in \textit{most} fundamental domains). Indeed, after applying the $z\mapsto -\frac{1}{z}$ transformation we get
    \[ g(z) \coloneqq f\left( -\frac{1}{z}\right) = \frac{4j(z)^2(j(z)-1728)^2 + \alpha z^2 j(z)^2(j(z)-1728)j'(z)}{z^6j'(z)^3}\]
    and because of the extra factor $z^2$ in the second summand in the numerator, no cancellation is possible, thus guaranteeing that $g(z)$ grows exponentially as $z \to i\infty$. Note however that this function is not periodic, but only \textit{asymptotically periodic} in the sense that $\lim_{z\to i\infty}\frac{g(z+1)}{g(z)}=1$. This fact is responsible for the technicalities in the rest of this section.

    It is worth mentioning that for equations of the form $F(j,j')=0$, where $F$ is a polynomial, the transformation $z\mapsto -\frac{1}{z}$ always guarantees that the ratio of the lowest and highest powers of $z$ has a pole in $\h$ or exponential growth at $i\infty$. In general, \Cref{prop: intro periodic pole} is sufficient when we deal with equations of the form $F(z,j(z),j'(z))=0$, although the argument is somewhat more complicated.

    The reader may benefit from revisiting this example after reading the rest of the paper, as it will make the above-mentioned phenomena less obscure.
\end{example}

\begin{notation}
  Let $\mathcal{P}$ denote the field of $1$-periodic meromorphic functions on $\h$ {which are also meromorphic at $i\infty$ (recall Definition~\ref{def: merom at infty})}. 
  We write $\mathcal{P}[w]$ and $\mathcal{P}(w)$ respectively for the polynomial ring and its fraction field generated by the variable $w$ over $\mathcal{P}$. 
  {We remark that the functions in $\mathcal{P}$ will also be thought of as meromorphic functions of the variable $w$.}

  {Also, given an unbounded region $U\subseteq\C$ and two meromorphic functions $f,g$ on $U$, we write \emph{$f\sim g$ for $w \to \infty$ in $U$} to mean that the limit of the ratio $\frac{f(w)}{g(w)}$ tends to 1 as $w$ approaches infinity from within $U$.}
\end{notation}

\begin{lemma}\label{lemma:exp-order}
    Let $f \in \mathcal{P}(w)$. Then there are $\alpha \in \C^\times$, $e, d \in \Z$, and a positive $C \in \R$ such that
    \[ f(w) \sim \alpha w^dq^e \]
    for $w \to \infty$ in the region $\im(w) \geq C\log|w|$, where $q = \exp(2\pi i w)$.
\end{lemma}
\begin{proof}
    It suffices to prove the conclusion for $f \in \mathcal{P}[w]$. Write $f(w) = \sum_{k=0}^n g_k(w)w^k$, with $g_k \in \mathcal{P}$.

    Each $g_k(w)$ has a meromorphic $q$-expansion $\tilde{g}_k(q)$ {which converges} on some neighbourhood of $q = 0$. Let $e$ be the minimum order of $\tilde{g}_k(q)$ at $q = 0$ for $k = 0, \dots, n$. {Let $\alpha_k \in \C$ be such that $\tilde{g}_k(q) = q^e(\alpha_k + O(q))$. Let $d$ be the maximum $k$ such that $\tilde{g}_k(q)$ has order $e$ at $q = 0$, namely such that $\alpha_k \neq 0$. In the region $|q| \leq |w|^{-(n+1)}$ we have $O(q) = O(w^{-(n+1)})$, in which case
    \[ G(w,q) = \sum_{k = 0}^n \tilde{g}_k(q)w^k = q^e w^d \sum_{k = 0}^n \left(\alpha_k + O\left(w^{-(n+1)}\right)\right)w^{k-d} = \alpha_d q^e w^d \left(1 + O\left(w^{-1}\right)\right). \]
    }

    It now suffices to specialise at $q = \exp(2\pi i w)$ and observe that $|q| = e^{-2\pi\im w} \leq |w|^{-(n+1)}$ if and only if $\im(w) \geq \frac{n+1}{2\pi}\log|w|$.
\end{proof}

It follows at once that all functions in $\mathcal{P}(w)$ are `asymptotically periodic' in the sense that $f(w+1) \sim f(w)$ for $w \to \infty$ in the above region. Moreover, the lemma allows us to make the following definition.

\begin{definition}
\label{def:expgrowth}
    Call the \emph{order at $i\infty$} of $f \in \mathcal{P}(w)$, written $\ord_{w=i\infty}(f)$, the pair $(e,d) \in \Z^2$ of integers such that for some $\alpha\in\C$ and some $C \in \R$, we have $f(w) \sim \alpha w^{-d}\exp(2\pi iew)$ for $w \to \infty$ in the region $\im(w) \geq C\log|w|$.

    We say that $f$ has \emph{exponential growth at $i\infty$} if its order is $(e,d)$ with $e < 0$.
\end{definition}
Here we consider $\Z^2$ as an ordered group, ordered lexicographically. This makes $(\mathcal{P}(w), \ord_{w=i\infty})$ into a valued field, and we have for instance that $f(w) \to \infty$ in a suitable region $\im(w) \geq C\log|w|$ if and only if $\ord_{w=i\infty} f(w) < (0,0)$.

{We can now set up a generalisation of \Cref{prop: periodic finite pole} that will cover our application to $j$. Let us fix the following data:
\begin{itemize}
    \item a polynomial $F(z,w) = \sum_{k=0}^n z^k f_k(w)$ where each $f_k$ is in $\mathcal{P}(w)$ and $f_n \neq 0$;
    \item a value $s$ which is either $0$ or $1$.
\end{itemize}}

\begin{figure}
    \centering
    \includegraphics{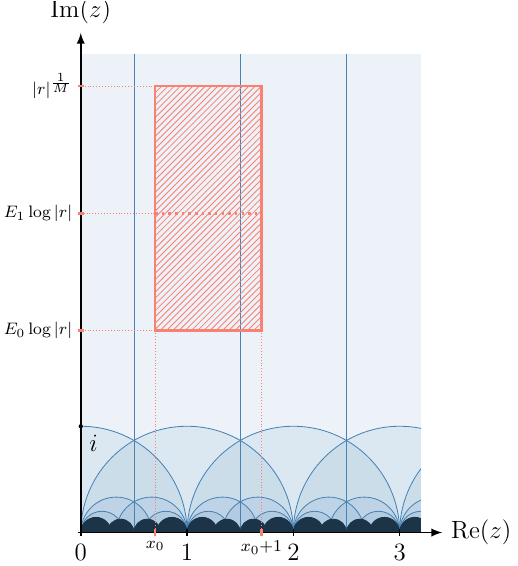}
    \caption{The region $\Xi_r$ highlighted by the striped background.}
    \label{fig:region-xi}
\end{figure}

{We look for the zeroes of functions of the form $F_r(w) \coloneqq F(r + sw, w)$ for $r$ varying among the real numbers. For each $r$ sufficiently large, we pick a suitable rectangle $\Xi_r$, pictured in \Cref{fig:region-xi} and defined in \Cref{zeroes if growth at infinity}, and integrate the logarithmic derivative $\frac{F_r'}{F_r}$ along the boundary of $\Xi_r$. Provided that $F_r$ does neither have zeroes nor poles on such a boundary, by the Argument Principle (see \Cref{thm: arg princ}), the value of the integral counts the difference between the number of zeroes and poles, with multiplicity, inside $\Xi_r$. We will choose $\Xi_r$ so that the integral has positive value.}

{We first parametrise the roots of $F(z,w)$ as a polynomial in $z$ in terms of $w$ varying in a suitable region. The resulting functions, which by construction are algebraic over $\mathcal{P}(w)$, admit an order at $i\infty$ which may be a pair of rational numbers, rather than only integers.}

\begin{lemma}\label{lem:algebraic-functions}
    {There is a positive $C \in \R$ such that in the region \[ U = \{ w \in \h :  \im(w) \geq C\log|w| \},\] there are holomorphic functions $\beta_1,\dots,\beta_n : U \to \C$ such that $F(\beta_k(w),w) = 0$ for all $w \in U$, and if $F(\beta,w) = 0$, then $\beta = \beta_k(w)$ for some $k$.}

    {Moreover, there are $\alpha_k \in \C^\times$, $e_k, d_k \in \Q$ such that $\beta_k(w) \sim \alpha_k w^{-d_k} q^{e_k}$ for $w \to \infty$ in $U$, where $q = \exp(2\pi i w)$ and $w^{-d_k} = \exp(-2\pi i d_k \log(w))$ for some holomorphic branch of $\log(w)$ on $U$.}
\end{lemma}
\begin{proof}
{    It suffices to prove the conclusion for $F$ irreducible as a polynomial over $\mathcal{P}(w)$.

    Let $F^{(z)} = \frac{\partial F}{\partial z}$. Then there are $G,H \in \mathcal{P}(w)[z]$ such that $GF + HF^{(z)} = 1$. We take $C$ large enough that by \Cref{lemma:exp-order}, the coefficients of $F$, $F^{(z)}$, $G$, and $H$ are holomorphic in the region $U = \{ w \in \h :  \im(w) \geq C\log|w| \}$. In particular, $F(z,w_0)$ and $F^{(z)}(z,w_0)$ have no common roots for any $w_0 \in U$, and so by the implicit function theorem, and because $U$ is simply connected, there are $m$ holomorphic functions $\beta_1,\ldots,\beta_m : U \to \C$ such that $F(\beta_t(w),w) = 0$ for every $t$ and $w \in U$, and moreover taking distinct values at all $w \in U$, thus if $F(\beta,w) = 0$, then $\beta = \beta_k(w)$ for some $k$.

    Now fix some $k$. For every $w \in U$, there is some $t$ such that $|\beta_k(w)^t f_t(w)| \geq |\beta_k(w)^\ell f_\ell(w)|$ for every $\ell$, and since $F(\beta_k(w),w) = 0$, there is also $h \neq t$ such that $|\beta_k(w)^h f_h(w)| \geq \frac{1}{n} |\beta_k(w)^t f_t(w)|$. Let $U_{t,h}$ be the region where those inequalities hold. We have in particular
    \[ 1 \geq |\beta_k(w)| \sqrt[h-t]{\left|\frac{f_h(w)}{f_t(w)}\right|} \geq \sqrt[h-t]{\frac{1}{n}} \]
    Let $(e_k,d_k) = \frac{\ord_{w=i\infty}(f_t / f_h)}{h - t}$. By \Cref{lemma:exp-order} combined with the above inequalities, there is $N > 1$ such that whenever $w$ is sufficiently large in $U_{t,h}$ we have $N \geq \frac{|\beta_k(w)|}{|w^{-d_k}q^{e_k}|} \geq \frac{1}{N}$ for some fixed determination of $w^{-d_k}$ on $U$. Choose $N$ so that it works simultaneously for any possible pair $t$, $h$. By continuity of $\beta_k$, the numbers $d_k$, $e_k$ do not depend on $t$, $h$.

    Let $S$ be the set of indices $t$ such that $\ord_{w=i\infty}(f_t) + (te_k,td_k)$ reaches a minimum value $(e,d)$. By construction, $S$ contains at least two elements. Write
    \[ F(w^{-d_k}q^{e_k}z',w) = q^ew^d(G_0(z') + G_1(z',w)) \]
    where now $G_0$ is a non-trivial polynomial in $z'$ with $|S| \geq 2$ terms and constant coefficients, and $G_1$ has coefficients that tend to $0$ for $w \to \infty$ in $U_{t,h}$. Since $\frac{\beta_k(w)}{w^{-d_k}q^{e_k}}$ is bounded and continuous, it must converge to a non-zero root $\alpha_k$ of $G_0(z')$, thus $\beta_k(w) \sim \alpha_k w^{-d_k} q^{e_k}$, as desired.}
\end{proof}

We now provide an estimate on the size of $F_r(w)$ that we can use on the boundary of $\Xi_r$.

\begin{lemma}\label{lem:lower bound F}
    There exist $0 \leq x_0 < 1$, $C, y_0 \geq 4$, $D, E_0, E_1 \geq 0$ (with $E_0$ possibly $+\infty$), $\delta > 0$ such that
    \[ |F_r(w)| = |F(r + sw, w)| > \delta\sum_{k=0}^{n-1}|f_k(w)||r + sw|^k \]
    for all $w \in \h$, $r \in \R$ such that $\im(w) \geq y_0$, $|r| \geq \im(w)^D$, and one of the following holds:
    \begin{itemize}
        \item $0 \leq \re(w) \leq 2$, $\im(w) \leq E_0\log|r|$, or
        \item $0 \leq \re(w) \leq 2$, $E_1\log|r| \leq \im(w)$, or
        \item $\re(w) \in \{x_0, x_0+1\}$.
    \end{itemize}
\end{lemma}
\begin{proof}
    {We work in a region $U = \{ w \in \h : 0 \leq \re(w) \leq 2 \wedge \im(w) \geq y_0 \}$ for some $y_0$ sufficiently large as determined by this proof. We start by taking $y_0$ large enough that by \Cref{lemma:exp-order},    $f_n$ has neither zeroes nor poles at $w$. We also require that $|r| \geq C$ for some $C$ sufficiently large, again as determined by this proof.

    By \Cref{lem:algebraic-functions}, provided $y_0$ is sufficiently large, there are holomorphic functions $\beta_1,\ldots,\beta_n : U \to \C$ parametrising the roots of $F(z,w) = 0$ as functions of $w$, and for $w \to \infty$ in $U$ we have
    \[ \beta_k \sim \alpha_k w^{-d_k} q^{e_k} \]
    for some $\alpha_k \in \C^\times$ and $d_k,e_k \in \Q$.

    Since we are assuming $0 \leq \re(w) \leq 2$, we also have $\im(w) \leq |w| \leq \im(w) + 2$. We require that $C \geq 4$, $y_0 \geq 4$, so that we have the following simple inequalities:
    \[ |r + sw| \geq \max\{|r| - 2, s\im(w)\} \geq \max\left\{\frac{|r|}{2},s\frac{|w|}{2}\right\} \geq \frac{|r + sw|}{4}. \]
    It follows, for instance, that $|w| \sim \im(w) \to +\infty$ for $w \to \infty$ in $U$. We shall omit the specification `in $U$' in the rest of this proof.

    We shall now bound $|r + sw - \beta_k(w)|$, distinguishing multiple cases depending on $(e_k,d_k)$.

    \begin{itemize}
        \item If $(e_k,d_k) \geq (0,0)$, then $\beta_k(w)$ converges to a finite value $\gamma_k$ for $w \to \infty$. In this case, we require $C$ (if $s = 0$) and $y_0$ to be large enough that
        \[ |r + sw| \geq 4|\gamma_k| \geq 2|\beta_k(w)|. \]
        This ensures that
        \[ |r + sw - \beta_k(w)| \geq \frac{|r + sw|}{2} \geq \frac{|r + sw| + |\beta_k(w)|}{4}. \]
        \item Suppose that $e_k = 0$ and $d_k < 0$. For $|r| \geq \im(w)^{2\max\{1,-d_k\}}$, we have both $|r| \geq 4|sw|$ and $|r| \geq 4|\beta_k(w)|$ for $w$ sufficiently large, and so for $y_0$ large we get
        \[ |r + sw - \beta_k(w)| \geq \frac{|r|}{2} \geq \frac{|r + sw| + |\beta_k(w)|}{8}. \]
     
        \item If $e_k < 0$, then $sw - \beta_k(w) \sim -\beta_k(w)$ for $w \to \infty$. Since $\re(w)$ is bounded, we get
        \[ \re(\log(sw - \beta_k(w))) = \log|sw - \beta_k(w)| \sim \log|\beta_k(w)| \sim -2\pi e_k\im(w). \]
        We first give bounds when $|r|$ is roughly at least $|\beta_k|^2$, and when $|\beta_k|$ is roughly at least $|r|^2$. More precisely, for $y_0$ sufficiently large, we have
        \begin{align*}
            \left|r + sw - \beta_k(w)\right| &> \frac{|r|}{2} \geq \frac{|r + sw| + |\beta_k(w)|}{8} &\text{for } |r| \geq e^{-4\pi e_k\im(w)}, \\
            \left|r + sw - \beta_k(w)\right| &> \frac{|\beta_k(w)|}{2} \geq \frac{|r + sw| + |\beta_k(w)|}{8} &\text{for } |r| \leq e^{-\pi e_k\im(w)}.
        \end{align*}
        For $w \to \infty$ we also have
        \[ \im(\log(sw - \beta_k(w))) \sim \arg(\alpha_k) + \pi + 2\pi e_k\re(w) \mod 2\pi.\]
       
        Now choose $x_0$ so that $\arg(\alpha_k) + 2\pi e_kx_0$ and $\arg(\alpha_k) + 2\pi e_k(x_0 + 1)$ are not in $\pi\Z$, and so are different from $\arg(r) \in \pi\Z$. Note that we can do this simultaneously for all $k$ such that $e_k < 0$. In particular, there is $\delta_k > 0$ such that, after taking $y_0$ sufficiently large, we have
        \[ \left|r + sw - \beta_k(w)\right| > \delta_k(|r + sw| + |\beta_k(w)|) \quad \text{for } \re(w) \in \{x_0,x_0+1\}. \]
    \end{itemize}
    Now, if $e_k = 0$ for some $k$, let $D$ be the maximum between the values $-2d_k$ and $2$ for such $k$, otherwise we can take $D = 0$. If $e_k < 0$ for some $k$, let $E_0$ be the maximum of $-4 \pi e_k$ and let $E_1$ be the minimum of $-\pi e_k$ for such $k$, otherwise let $E_0 = +\infty$ and $E_1 = 0$. Under the above choices, there is $\delta > 0$ such that
    \[ |F_r(w)| > \delta |f_n(w)| \prod_k (|r + sw| + |\beta_k(w)|) \geq \delta \sum_{k=0}^n |f_k(w)||r + sw|^k \]
    for any $w \in \h$, $r \in \R$ satisfying the requirements in the conclusion.}
\end{proof}

We now recall the Argument Principle, which plays a key role in the proof of \Cref{zeroes if growth at infinity}.

\begin{theorem}[Argument Principle, see e.g.~{{\cite[Chapter VI, \S1, Theorem 1.5]{lang:complexanalysis}}}]\label{thm: arg princ}
    {Let $f$ be a meromorphic function on a complex domain $\Omega$.
    Let $C$ be a simple closed curve (positively oriented) which is homologous to $0$ in $\Omega$ and such that $f$ has no zeroes or poles on $C$. Let $Z$ and $P$ respectively denote the number of zeroes and poles (counted with multiplicity) of $f$ in the interior of $C$. Then
    \[ 2\pi i (Z-P) = \oint_C \frac{f'(z)}{f(z)}dz = \oint_{f \circ C} \frac{dz}{z}. \]}
\end{theorem}

\begin{figure}
    \centering
    \includegraphics{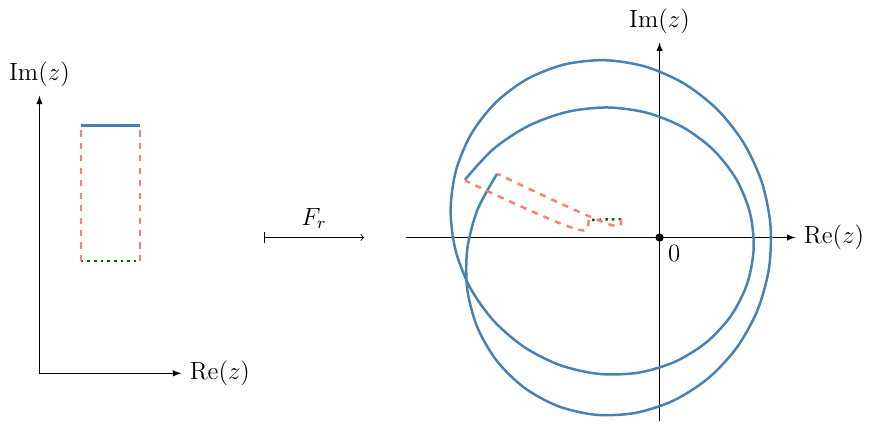}
    \caption{Visual representation of the action of $F_r$ on the boundary of a typical rectangle $\Xi_r$ for $F$ of the form $Az^2 + Bz + C{j(z)}^2 + Dj(z) + E$. The term $j^2$ has lowest order $(-2,0)$, so $F_r$ winds around $0$ twice while following the top side of the rectangle.}
    \label{winding}
\end{figure}

{In the proof of the following proposition, we will integrate $\frac{F_r'}{F_r}$ along the boundary of $\Xi_r$ and use the above estimates to find a positive lower bound, thus proving the existence of zeroes of $F_r$ within $\Xi_r$. For a more geometric description, integrating $\frac{F_r'}{F_r}$ computes how many times the image $F_r(z)$ winds around $0$ while $z$ moves along $\partial \Xi_r$. The bounds below will determine a rough picture of $F_r(\partial \Xi_r)$, as in \Cref{winding}, and in turn determine the number of zeroes, counted with multiplicity.}

\begin{proposition}\label{zeroes if growth at infinity}
    Let $(e,d)$ be the minimum order of $\frac{f_k}{f_n}$ at $i\infty$ for $k = 0, \dots, n - 1$ and suppose that $e < 0$. We work under the notation of \Cref{lem:lower bound F}.

    Then for all $r \in \R$ sufficiently large the function $F(r + sw,w)$ has $-e$ zeroes, counted with multiplicity, within the region (see \Cref{fig:region-xi})
    \[ \Xi_r = \{ w \in \h : x_0 < \re(w) < x_0 + 1,\ E_0\log|r| < \im(w) < |r|^{\frac{1}{M}} \}, \]
    where $M = \lceil D \rceil$ if $D > 0$ and $M = 1$ otherwise. Moreover, for $w \in \partial\Xi_r$ we have
    \[ |F(r + sw, w)| \geq \delta \sum_{k=0}^n |f_k(w)||r + sw|^k. \]
\end{proposition}

\begin{proof}
    Recall that $F_r(w) = F(r + sw, w)$. First, we apply \Cref{lem:lower bound F} to $F_r(w)$ and find relevant constants $x_0$, $y_0$, $\delta$, $C$, $D$, $E_0$, $E_1$. Let $M = \lceil D \rceil$ if $D > 0$ and $M = 1$ otherwise. Let $r$ be large enough so that $|r| \geq C$, $E_0\log|r| \geq y_0$, and $E_1\log|r| \leq |r|^{\frac{1}{M}}$, hence
    \[ |F_r(w)| \geq \delta\sum_{k=0}^n|f_k(w)||r + sw|^k \]
    for $w \in \partial\Xi_r$. It follows, for instance, that $F_r(w)$ has order $(e,d)$ at $i\infty$. We also take $r$ sufficiently large so that each $f_k(w)$ does not have poles in $\Xi_r$, so in particular $F_r(w)$ is holomorphic on $\Xi_r$.

    Note that $F_r(w)$ is never zero on the boundary of $\Xi_r$, thus its logarithmic derivative $\frac{F_r'}{F_r}$ is holomorphic there. We shall now compute the integral of $\frac{F_r'}{F_r}$ of $F_r$ along such a boundary.

    \smallskip\noindent\emph{Vertical sides.} {We show that the images of the vertical sides of $\Xi_r$ under $F_r$ must be close to each other, and so their contributions cancel out as they are taken with opposite orientations. For instance, if $F_r$ happens to be $1$-periodic (as a function of $w$; for instance, when $s = 0$ and the coefficients of $F$ are $1$-periodic) then the images of these vertical sides are identical.}

    First, we observe that
    \[ F_r(w+1) - F_r(w) = \sum_{k=0}^n f_k(w)(r + sw)^k\left(\frac{f_k(w+1)}{f_k(w)}\left(\frac{r + sw + 1}{r + sw}\right)^k - 1\right). \]
    Since $f_k \in \mathcal{P}(w)$, we have that $f_k(w + 1) \sim f_k(w)$ as $\im(w) \to +\infty$. Likewise, $r + sw + 1 \sim r + sw$ for $r + sw \to \infty$. Therefore, we can choose $r$ large enough so that the last factor on the right-hand side has modulus less than $\frac{\delta}{2}$ for all $k$ and for any $w$ on the boundary of $\Xi$. We then have
    \begin{equation}\label{Fr-small-changes}
        \left|\frac{F_r(w+1)}{F_r(w)} - 1\right| = \left|\frac{F_r(w+1) - F_r(w)}{F_r(w)}\right| \leq \frac{\delta}{2|F_r(w)|}\left(\sum_{k=0}^n|f_k(w)||r + sw|^k\right) < \frac{1}{2}
    \end{equation}
    for $w \in \partial\Xi$. We may now choose a branch of $\log$ in the disc around $1$ of radius $\frac{1}{2}$ and estimate the integral along the vertical sides as\footnote{{Recall that $\frac{(fg)'}{fg} = \frac{f'}{f} + \frac{g'}{g}$, and that $\int_\gamma \frac{f'}{f} = \log(f(\gamma(1))) - \log(f(\gamma(0)))$ whenever $\gamma : [0,1] \to \C$ is a path and $f \circ \gamma$ takes values in a simply connected region on which we have fixed a continuous branch of $\log$}.}
    \begin{multline*}
        \left|\int_{E_0\log|r|}^{|r|^{\frac{1}{M}}} \left(\frac{F_r'(x_0 + 1 + iy)}{F_r(x_0 + 1 + iy)} - \frac{F_r'(x_0 + iy)}{F_r(x_0 + iy)}\right) dy\right| = \\
        = \left|\left.\log\left(\frac{F_r(w+1)}{F_r(w)}\right)\right|_{w = x_0 + iE_0\log|r|}^{w = x_0 + i|r|^{\frac{1}{M}}}\right| < \log\left(\frac{3}{2}\right) - \log\left(\frac{1}{2}\right) = \log(3) < 2.
    \end{multline*}

    \smallskip\noindent\emph{Bottom side.} {We now show that the image of the bottom side is away from $0$ and cannot wind much.} Using \eqref{Fr-small-changes}, for $r$ sufficiently large,
    \[ \left|\int_{x_0}^{x_0+1} \frac{F_r'(x + iE_0\log|r|)}{F_r(x + iE_0\log|r|)} dx\right| = \left|\log\left(\frac{F_r(x+iE_0\log|r|+1)}{F_r(x+iE_0\log|r|)}\right)\right| < \log\left(\frac{3}{2}\right) < 1. \]

    \smallskip\noindent\emph{Top side.} {We show that $F_r$ behaves like $\exp(2\pi i e w)$ on the top side of $\Xi_r$, and so the image under $F_r$ is roughly a circle traversed approximately $e$ times.}

    We now constraint $w$ to the region $\im(w) = |r|^{\frac{1}{M}}$, $x_0 \leq \re(w) \leq x_0 + 1$. We have
    \[ \delta\max_k|f_k(w)||r + sw|^k \leq |F_r(w)| \leq n\max_k|f_k(w)||r + sw|^k. \]
    By construction, $r \sim \zeta w^M$ for some power $\zeta$ of $i$ depending on $M$ and the sign of $r$. For simplicity, fix the sign of $r$, and assume that $\zeta = 1$, so that $r \sim w^M$. Then
    \[ F(w^M + sw, w) - F_r(w) = \sum_{k=0}^n f_k(w)(r + sw)^k\left(\left(\frac{w^M + sw}{r + sw}\right)^k - 1\right). \]
    In particular, by \Cref{lem:lower bound F}, we find that $F(w^M + sw, w) - F_r(w) = o(F_r(w))$ for $r \to +\infty$. Since $F(w^M + sw, w)$ is in $\mathcal{P}(w)$, it has an order $(e',d')$ at $i\infty$, and in fact $e' = e$ because the term $w^M$ cannot alter the exponential growth.

    Therefore, we find that there is $\alpha \in \C^{\times}$ such that for $r$ large enough
    \[ \left|\frac{F_r(w)}{\alpha w^{-d'}\exp(2\pi iew)} - 1 \right| < \frac{1}{4}. \]

    Observe that for $r$ large enough we have
    \[ \left|\int_{x_0}^{x_0 + 1}\frac{((x + i|r|^{\frac{1}{M}})^{-d'}\exp(2\pi ie(x + i|r|^{\frac{1}{M}})))'}{(x + i|r|^{\frac{1}{M}})^{-d'}\exp(2\pi ie(x + i|r|^{\frac{1}{M}}))}dx - 2\pi ie\right| = \left|\int_{x_0}^{x_0 + 1}\frac{d'}{x + i|r|^{\frac{1}{M}}}dx\right| < 1. \]
    Thus, for sufficiently large $r$ we have
    \[ \left|\int_{x_0}^{x_0 + 1} \frac{F_r'(x + i|r|^{\frac{1}{N}})}{F_r(x + i|r|^{\frac{1}{N}})} dx - 2\pi ie\right| \leq \left|\left.\log\left(\frac{F_r(w)}{\alpha w^{-d'}\exp(2\pi iew)}\right)\right|_{w = x_0 + i|r|^{\frac{1}{N}}}^{w = x_0 + 1 + i|r|^{\frac{1}{N}}}\right| + 1 < 2. \]

    \smallskip\noindent\emph{Conclusion.}   {Using the above estimates we can now find the winding number of $F_r(\partial \Xi_r)$ at $0$.}

    Summing up the contributions from all sides, with the appropriate orientations, we obtain
    \[ \left|\frac{1}{2\pi i}\oint_{\partial \Xi_r} \frac{F_r(w)'}{F_r(w)} dw + e\right| < \frac{1}{2\pi}\left(2 + 1 + 2\right) < 1. \]

    By the Argument Principle, the integral on the left-hand side must be the difference between the number of zeroes and poles of $F_r(w)$ inside $\Xi$ (in particular an integer) multiplied by $2\pi i$. Since $F_r(w)$ is holomorphic on $\Xi_r$, it must have $-e$ zeroes in $\Xi_r$, counted with multiplicity.
\end{proof}

\begin{proof}[Proof of Proposition \ref{prop: intro periodic pole}]
    This follows from combining Propositions \ref{prop: periodic finite pole} and \ref{zeroes if growth at infinity}, where we use $r=m$ a large integer and $s=1$. Note that if $\frac{f_k}{f_n}(z) \to \infty$ as $z\to i\infty$ then $\frac{f_k}{f_n}(z)$ must have exponential growth at $i\infty$ for it is periodic and so has a $q$-expansion.
\end{proof}

\section{Some Criteria for Zariski Density}
\label{sec:zariskidensity}

 Recall that  $\Y:=(Y_0,Y_1,Y_2)$ and $\J:=(j,j',j''): \h \to \C^3$.

\subsection{Generic transforms}
\label{subsec:gentransform}
{Given $p \in \C[X,\Y]$, or more generally $p \in K[X,\Y]$ for some field $K$, we define the \textit{generic $\SL_2(\Z)$-transform of $p$} to be the polynomial $\Gamma(p) \in K[Z,W,C,\mathbf{Y}]$ given by
\[ \Gamma(p)(Z,W,C,\Y) \coloneqq p\left(Z,Y_0,W^2 Y_1,W^4 Y_2 + 2CW^3 Y_1\right). \]
In particular $\deg_X(p) = \deg_Z(\Gamma(p))$.
By construction, for any $\gamma = \begin{pmatrix} a & b \\ c & d \end{pmatrix} \in \SL_2(\Z)$, if $n = \deg_X(p)$ we have
\[ p(\gamma z,\J(\gamma z)) = \frac{p^\gamma(z,\J(z))}{(cz+d)^n}, \ \text{where }p^{\gamma}(X,\Y) \coloneqq (cX+d)^n \Gamma(p)\left(\frac{aX+b}{cX+d},cX+d,c,\Y\right). \]
Note that $p^\gamma \in K[X,\Y]$.}

{We make the following observations.
\begin{enumerate}[label={\bf O\arabic*.},ref={O\arabic*}]
    \item\label{Gamma-homomorphism} The map $\Gamma:K[X,\Y]\to K[Z,W,C,\Y]$ defined above is a $K$-algebra homomorphism with left inverse $p(X,\Y) = \Gamma(p)(X,1,0,\Y)$.

    \item\label{gamma-multiplicative} For each $\gamma\in\SL_2(\Z)$, the map $p\mapsto p^\gamma$ is multiplicative, that is $(p_1 p_2)^\gamma = p_1^\gamma p_2^\gamma$ for any $p_1,p_2\in K[X,\Y]$.    Indeed, this follows at once from the fact that $\Gamma$ is a homomorphism \eqref{Gamma-homomorphism} and $\deg_X(p_1p_2) = \deg_X(p_1) + \deg_X(p_2)$.

    \item\label{gamma-homomorphism-up-to-r} For any $p\in K[X,\Y]$ and any $\gamma_1,\gamma_2 \in \SL_2(\Z)$ there is $r \in \Z[X]$ such that $(p^{\gamma_1})^{\gamma_2} = r(X)p^{\gamma_1\gamma_2}$.

    Indeed, let $n = \deg_X(p)$ and $m = \deg_X(p^{\gamma_1})$. By construction, we have $m = n + \deg_{Y_0}(p) \geq n$. Now write $\gamma_t = \begin{pmatrix} a_t & b_t \\ c_t & d_t \end{pmatrix}$, for $t\in\{1,2\}$, and write $\gamma_1\gamma_2 = \begin{pmatrix} \widetilde{a} & \widetilde{b} \\ \widetilde{c} & \widetilde{d} \end{pmatrix}$.
    Thus
    \[\frac{p^{\gamma_1\gamma_2}(z,\J(z))}{\left(\widetilde{c}z+\widetilde{d}\right)^n} = p(\gamma_1\gamma_2z,\J(\gamma_1\gamma_2z)) = \frac{p^{\gamma_1}(\gamma_2z,\J(\gamma_2z))}{(c_1\gamma_2z+d_1)^n} = \frac{\left(p^{\gamma_1}\right)^{\gamma_2}(z,\J(z))}{(c_1\gamma_2z+d_1)^n(c_2z+d_2)^m}.\]
    Since $m\geq n$ and $(c_1\gamma_2z+d_1)(c_2z+d_2) = \widetilde{c}z+\tilde{d}$, we get
    \[r(X):=\frac{(c_2X+d_2)^m}{\left(\widetilde{c}X+\widetilde{d}\right)^n}\left(c_1\frac{a_2X+b_2}{c_2X+d_2} + d_2\right)^n = (c_2X+d_2)^{m-n}\in\Z[X].\]

    \item\label{gamma-irreducibility} For any $p\in K[X,\Y]$ and any $\gamma\in\SL_2(\Z)$, if we consider $p$ and $p^\gamma$ as polynomials in the variables $\Y$ with coefficients in $K(X)$, then $p$ is irreducible in $K(X)[\Y]$ if and only if so is $p^\gamma$.

    Indeed, note that if $p$ is not a unit (meaning it contains one of the variables $Y_0$, $Y_1$, $Y_2$), then $p^\gamma$ is also not a unit. It follows by \eqref{gamma-multiplicative} that if $p$ is reducible in $K(X)[\Y]$, thus it is a product of two non-units, then so is $p^\gamma$. Likewise, if $p^\gamma$ is reducible, then $(p^\gamma)^{\gamma^{-1}}$ is reducible too, and $(p^\gamma)^{\gamma^{-1}} = r(X)p$ for some $r(X) \in \C[X]$ by \eqref{gamma-homomorphism-up-to-r}; since $r(X)$ is a unit, it follows that $p$ is reducible.
\end{enumerate}}

\begin{proposition}
\label{prop:p^gamma}
    Let $p$ be an irreducible polynomial in $\C[X,\Y] \setminus \C[X]$ and $\gamma \in \SL_2(\Z)$. Let $h$ be the irreducible factor of $p^\gamma$ that is not in $\C[X]$. Then the equation $p(z,\J(z)) = 0$ has a Zariski dense set of solutions if and only if the equation $h(z,\J(z)) = 0$ has a Zariski dense set of solutions.
\end{proposition}
\begin{proof}

    {By \eqref{gamma-homomorphism-up-to-r} and \eqref{gamma-irreducibility}, $p^\gamma = r(X)h$ and $h^{\gamma^{-1}} = s(X)p$ for some $r, s \in \C[X]$. It follows that $p(z,\J(z)) = 0$ and $h(z,\J(z)) = 0$ have the same solutions except possibly for the zeroes of $r(z)$ and $s(z)$. Since those are only finitely many, the solutions of the former equation are Zariski dense if and only if so are the solutions of the latter.}
\end{proof}

\subsection{Density criteria}
\label{subsec:densitycriteria}
We now apply the results of \S\ref{sec:periodic} to establish some useful criteria for Zariski density of solutions of equations involving $z,j(z),j'(z),j''(z)$.

\begin{definition}
    Given a function $g\in\C(z,\J(z))$ and $\gamma\in\SL_2(\Z)$, we say that \emph{$g(z)$ has exponential growth in $\gamma\F$}, if $g\left(\gamma^{-1}z\right)$ has exponential growth at $i\infty$. Furthermore, if $r$ is the cusp of $\gamma\F$ (that is, $r \in \Q \cup \{\infty\}$ is in the Euclidean closure of $\gamma\F$), then we define the \emph{order of $g(z)$ in $\gamma\F$ at $r$} as $\ord_{z=i\infty}\left(g\left(\gamma^{-1}z\right)\right)$.
\end{definition}

\begin{proposition}\label{exp growth implies solutions}
    Let $F(X,\Y) = \sum_{k = 0}^n X^kp_k(\Y)$ be a polynomial. Assume that for some $k$ and some $\gamma \in \SL_2(\Z)$, the function $\frac{p_k}{p_n}(\J(z))$ has exponential growth in $\gamma \F$.

    Then there are $\ell > 0$, $0 \leq x_0 < 1$,
    $M > 0$, $E_0 > 0$ such that for all $m \in \Z$ sufficiently large the function $F(z,\J(z))$ has $\ell$ zeroes, counted with multiplicity, within the region $m + \gamma\Xi_m$, where
    \[ \Xi_m = \{ z \in \h : x_0 < \re(z) < x_0 + 1,\ E_0\log|m| < \im(z) < |m|^{\frac{1}{M}} \}. \]
\end{proposition}
\begin{proof}
    Fix some $s \in \{0,1\}$, $t \in \R$ to be determined later. For $m \in \Z$, let
    \[ F_m(z) \coloneqq F(m + t + sz, \J(\gamma z)). \]
    Likewise, set $G(z) \coloneqq F(z,\J(z))$. By \Cref{zeroes if growth at infinity} applied to $F(z,\J(\gamma z))$, $F_m(z)$ has $\ell > 0$ zeroes in a certain region $\Xi_m$ and is suitably bounded from below for $z \in \partial\Xi_m$, as long as $m$ is sufficiently large.

    If $\gamma$ is upper triangular, namely $\gamma z = z + k$ for some $k$, we choose $s = 1$, $t = 0$, and observe that since the functions of $\J$ are 1-periodic
    \[ G(z + m) = F(z + m,\J(\gamma z)) = F_m(z), \]
    thus $G(z)$ has $\ell$ zeroes in each region $m + \Xi_m$.

    Otherwise, let $s = 0$ and let $t$ be the limit of $\gamma z$ as $z \to \infty$ (where in fact $t \in \Q$). In particular,
    \begin{align*}
        G(m + \gamma z) - F_m(z) &= F(m + \gamma z,\J(\gamma z)) - F(m + t,\J(\gamma z)) \\
        &= \sum_{k=0}^n p_k(\J(\gamma z))(m + t)^k\left(\left(\frac{m + \gamma z}{m + t}\right)^k - 1\right).
    \end{align*}
    Thus, as soon as $z$ is sufficiently large, the last factor on the right-hand side has modulus less than $\frac{1}{2}$ independently of $m$. Then pick $m$ large enough so that this happens whenever $\im(z) > E_0\log|m|$, and so
    \[ |G(m + \gamma z) - F_m(z)| < \frac{1}{2}|F_m(z)|. \]
    Therefore, by Rouché's theorem \ref{thm:rouche}, $G(m + \gamma z)$ and $F_m(z)$ have the same number of zeroes in $\Xi_m$, counted with multiplicity. It follows that $G(z)$ has $\ell$ zeroes in the region $m + \gamma\Xi_m$.
\end{proof}

\begin{proposition}\label{cor: pole at cusp Zariski density}\label{cor:pole j,j',j'' - dense}
    Let $F(X,\Y) = \sum_{k = 0}^n X^kp_k(\Y)$ be irreducible. Suppose that for some $k$, the function $P(z) = \frac{p_k(\J(z))}{p_n(\J(z))}$ satisfies one of the following:
    \begin{enumerate}[label=(\roman*)]
        \item $P(z)$ has a pole in $\h$, or
        \item $P(z)$ has exponential growth in some fundamental domain.
    \end{enumerate}
    Then the equation $F(z,\J(z)) = 0$ has a Zariski dense set of solutions.
\end{proposition}
\begin{proof}
    Suppose by contradiction that all the solutions of $F(z,\J(z)) = 0$ lie on a further hypersurface $G = 0$, where $G\in\C[X,\Y]$ is a non-constant polynomial not divisible by $F$. In particular, the algebraic subset of $\C^4$ defined by $\{F = G = 0\}$ has dimension two, so its projection onto the variables $Y_0,Y_1,Y_2$ has dimension at most two, meaning that the solutions satisfy an equation $H(\J(z)) = 0$ for some non-constant polynomial $H\in\C[\Y]$. The assumption on $F$ implies that $F$ depends on the variable $X$ (i.e.~$n\geq 1$), thus $F$ and $H$ are coprime.

    By Propositions~\ref{prop: periodic finite pole} and ~\ref{exp growth implies solutions}, for $m \in \Z$ large, there are regions $\Xi_m$ such that the original equation has solutions in $m + \Xi_m$, and moreover the real part of each $\Xi_m$ is bounded from above and below and the imaginary part is bounded away from $0$ uniformly in $m$. For each solution $\tau \in m + \Xi_m$, there are only finitely many integers $k$ such that $\tau + k$ is also a solution. This implies that for some $m$ sufficiently large, the function $H(\J(z))$ has infinitely many zeroes in the region $\bigcup_{|k| > |m|}\Xi_k$.
    If this union is bounded, then we conclude that $H$ is constantly zero by the identity theorem from complex analysis, but this contradicts the algebraic independence of $j$, $j'$, $j''$.
    So we assume that the union is unbounded, but by \Cref{lemma:exp-order} there exist $\alpha\in\C^\times$, $d,e\in\Z$ and $C>0$ such that $H(z)\sim \alpha z^d\exp(e2\pi iz)$ in the region $U:=\{z\in\h : \mathrm{Im}(z)\geq C\log|z|\}$.
    Then the only way $H$ can have infinitely many zeroes in $\bigcup_{|k| > |m|}\Xi_k$ is if $H$ is constantly zero, again contradicting the algebraic independence of $j$, $j'$, $j''$.
    This completes the proof.
\end{proof}

\begin{corollary}\label{cor:Zariski density begets Zariski density}
    Let $F(X,\Y) = \sum_{k = 0}^n X^kp_k(\Y)$ be irreducible. Suppose that $p_n$ has a factor $h$ such that the equation $h(\J(z)) = 0$ has a Zariski dense set of solutions. Then the equation $F(z,\J(z)) = 0$ has a Zariski dense set of solutions.
\end{corollary}
\begin{proof}
    {If $n=0$, then $F = p_0$ which, by irreducibility, is equal to a constant multiple of $h$, and so the result is immediate.
    If instead $n>0$, then by irreducibility of $F$ for some $k\in\{0,\ldots,n-1\}$, $p_k$ is non-zero and not divisible by $h$.
    Then $\frac{p_k(\J(z))}{p_n(\J(z))}$ will have poles in $\h$ at those solutions of $h(\J(z)) = 0$ which satisfy $p_k(\J(z))\neq 0$ (which exist since we are assuming Zariski density of the solutions of $h(\J(z)) = 0$).
    So now the corollary follows from \Cref{cor: pole at cusp Zariski density}.}
\end{proof}

\section{Zero Estimates}
\label{sec:zeroesetimates}

In view of the results in the previous section, we will now look at the poles of quotients of the form $\frac{p_k(\J)}{p_n(\J)}$, with $p_k,p_n\in\C[\Y]$. We keep using the notation introduced in \S\ref{subsec:gentransform}.

Given $p \in \C[X,\Y]$, we will prove a few estimates on the order of $p(z,\J(z))$ at different points, distinguishing three cases. Before that, we note that specialising the variables $Y_1$ and $Y_2$ of $\Gamma(p)$ at some complex values will almost always return an `obfuscated' copy of the original polynomial, which for instance cannot be constant unless $p$ itself was. More precisely, we note the following trivial identity.

\begin{lemma}\label{lem:specialise generic transform}
    For every $\alpha,\beta \in \C$ with $\alpha \neq 0$ we have
    \[ \Gamma(p)\left(X,\alpha^{-1} U_1,\alpha\frac{U_2 - \alpha^{-4}\beta U_1^4}{2U_1^3},Y_0,\alpha^2,\beta\right) = p(X,Y_0,U_1^2,U_2). \]
\end{lemma}
\begin{proof}
    Immediate.
\end{proof}

First, we consider unramified points of $j$, that is, points $\tau \in \h$ such that $j'(\tau) \neq 0$. If $(Y_0 - j(\tau))^s$ divides $p$, then $p(z,\J(z))$ has obviously order at least $s$ at all conjugates of $\tau$. We show that the order is exactly $s$ for most conjugates.

\begin{proposition}\label{prop:generic order at unramified tau}
    Let $p \in \C[X,\Y]$ be non-zero, $\tau \in \h$. Suppose that $j'(\tau) \neq 0$ and let $s$ be the maximum integer such that $(Y_0 - j(\tau))^s$ divides $p$. Then for all $\gamma$ in a Zariski open dense subset of $\SL_2(\Z)$, the function $p(z,\J(z))$ has order $s$ at $z = \gamma\tau$.
\end{proposition}
\begin{proof}
 {   Let $f\in\C[X,\Y]$ be such that $p = (Y_0-j(\tau))^sf$; we need to show that $f(z,\J(z))$ has order 0 at $z=\gamma\tau$, for all $\gamma$ in a Zariski open dense subset of $\SL_2(\Z)$.
    In other words, it suffices to prove the proposition for the case when $p$ is not divisible by $Y_0 - j(\tau)$ (and hence $s=0$).}

    From now on we assume $s=0$.
    By \Cref{lem:specialise generic transform} at $\alpha^2 = j'(\tau) \neq 0$, $\beta = j''(\tau)$, $U_1^2 = Y_1$, $U_2 = Y_2$, there are $V_1, V_2 \in \C(U_1,U_2)$ such that
    \[ \Gamma(p)\left(X,V_1,V_2,\J(\tau)\right) = p(X,j(\tau),Y_1,Y_2). \]
    Since $p$ is not divisible by $Y_0 - j(\tau)$, $p(X,j(\tau),Y_1,Y_2)$ is a non-zero polynomial, hence in particular $r(Z,W,C) \coloneqq \Gamma(p)(Z,W,C,\J(\tau))$ is also non-zero.

    Now take $\gamma = \begin{pmatrix} a & b \\ c & d \end{pmatrix}\in \SL_2(\Z)$ and write
    \[ p^\gamma(\tau,\J(\tau)) = (c\tau + d)^n r(\gamma\tau,c\tau + d,c) = (c\tau + d)^n p(\gamma\tau,\J(\gamma \tau)), \]
    where $n=\deg_X(p)$.
    In particular, $p(\gamma\tau,\J(\gamma\tau)) = 0$ if and only if $r(\gamma\tau,c\tau+d,c) = 0$. 
    
    {The map $\upsilon:\SL_2(\C)\to\C^3$ given by $\gamma \mapsto (\gamma\tau,c\tau + d,c)$ is injective, and since $\dim\SL_2(\C) = 3 = \dim\C^3$, by the fibre-dimension theorem $\upsilon$ is also dominant.
    As $r$ is a non-zero polynomial, there is a non-empty Zariski open subset $U_1$ of $\C^3$ such that $r$ never vanishes on $U$.
    This then gives a Zariski open subset $U_0$ of $\SL_2(\C)$ such that for every $\gamma\in U_0$, $\upsilon(\gamma)\in U_1$.
    This implies that for $\gamma\in U_0$, $p(z,\J(z))$ does not vanish at $z = \gamma\tau$, and hence $p(z,\J(z))$ has order $s = 0$ at $z=\gamma\tau$.}
\end{proof}

With this proposition we obtain the following special case of \Cref{thm: intro - main thm}.

\begin{corollary}\label{cor:j(z)-u Zariski dense}
    The equation $j(z) - u = 0$ has a Zariski dense set of solutions if and only if $u \notin \{0,1728\}$.
\end{corollary}
\begin{proof}
    {If $u \notin \{0,1728\}$, then for any $\tau\in\h$ satisfying $j(\tau)=u$ we have that $j'(\tau)\neq 0$.
    We need to show that the solutions of $j(z)=u$ are Zariski dense.
    So suppose that $p(X,\Y)\in\C[X,\Y]$ is such that its zero locus contains all the solutions of $j(z)=u$.
    The solutions of $j(z)=u$ are precisely $\SL_2(\Z)\tau$, where $j(\tau)=u$, so $p(\gamma\tau,\J(\gamma\tau))=0$ for all $\gamma\in\SL_2(\Z)$.
    Therefore $p(z,\J(z))$ has positive order at $\gamma\tau$ for all $\gamma\in\SL_2(\Z)$, so by \Cref{prop:generic order at unramified tau}, $Y_0-u$ divides $p$, thus proving Zariski density.

    On the other hand, if $u\in\{0,1728\}$, then for any $\tau\in\h$ satisfying $j(\tau)=u$ we have that $j'(\tau)= 0$, so the solutions of $j(z)=u$ lie in the proper Zariski closed subset given by $Y_1=0$.}
\end{proof}

Now we consider the behaviour of $p(z,\J(z))$ towards the cusps of the fundamental domains.
We write $T\Y:=(TY_0,TY_1,TY_2)$.
One can easily see that the order at the cusp of any given fundamental domain is at least $(-e,-N)$ for some $N$, where $e = \deg_T(p(X,T\Y))$. This is not far from the actual behaviour at most cusps.

\begin{proposition}\label{prop:generic order at infinity}
    Let $p \in \C[X,\Y]$. Then there is $0 \leq M \leq \deg_W(\Gamma(p))$ such that for all $\gamma$ in a Zariski open dense subset of $\SL_2(\Z)$, the function $p(z,\J(z))$ has order $(-e,-M)$ at the cusp of $\gamma\F$, where $e$ is the degree of $p(X,T\Y)$ in $T$.
\end{proposition}
\begin{proof}
    Let $n = \deg_X(p)$, $e = \deg_T(p(X,T\Y))$, and write
    \[ p(X,T\Y) = \sum_{k=0}^e T^k p_k(X,\Y) \]
    where each $p_k$ is homogeneous of degree $k$ in the variables $\Y$, namely
    \[ p_k(X,T\Y) = T^k p_k(X,\Y). \]
    Since $\Gamma(p)(Z,W,C,T\Y) = \Gamma(p(X,T\Y))$, we have
    \[ \Gamma(p)(Z,W,C,T\Y) = \sum_{k=0}^e T^k \Gamma(p_k)(Z,W,C,\Y), \]
    and so each $\Gamma(p_k)$ is still homogeneous of degree $k$ in $\Y$.

    Recall that for $\im(z) \to +\infty$, letting $q = \exp(2\pi iz)$ we have
    \[ j(z) = q^{-1} + O(1), \quad j'(z) = -2\pi i q^{-1} + O(1), \quad j''(z) = -4\pi^2 q^{-1} + O(1). \]
    Now fix some arbitrary $\gamma = \begin{pmatrix} a & b \\ c & d \end{pmatrix}\in \SL_2(\Z)$. We have
    \[ p^\gamma(z,\J(z)) = q^{-e} (cz+d)^n r(\gamma z,cz+d,c) + O(q^{-e+1}z^{n+N})  \]
    for $\im(z) \to +\infty$ in the standard fundamental domain, where $N$ is the degree of $\Gamma(p)$ in the variable $W$, and
    \[ r(Z,W,C) \coloneqq \Gamma(p_e)(Z,W,C,1,-2\pi i,-4\pi^2). \]
    More precisely, when $c \neq 0$, letting $M = \deg_W(r)$, we can also say
    \[ p^\gamma(z,\J(z)) = q^{-e}(cz)^n r(\sfrac{a}{c},cz,c) + O(q^{-e}z^{n+M-1}). \]
    Therefore, the order of $p^\gamma(z,\J(z))$ at $i\infty$ is $(-e,-(n+M))$, and so the order of $p(\gamma z,\J(\gamma z))$ is $(-e,-M)$, unless the leading coefficient of $r(\sfrac{a}{c},W,c)$ vanishes or $c = 0$.

    To conclude, since the map $\gamma \mapsto (a,c)$ from $\SL_2(\C)$ to $\C^2$ is dominant, it suffices to prove that $r$ is non-trivial (like we did in the proof of \Cref{prop:generic order at unramified tau}). By \Cref{lem:specialise generic transform} at $\alpha^2 = -2\pi i$, $\beta = -4\pi$, $U_1^2 = Y_1 Y_0^{-1}$, $U_2 = Y_2 Y_0^{-1}$, there are $V_1,V_2 \in \C(U_1,U_2)$ such that
    \begin{align*}
        r(X,V_1,V_2) = \Gamma(p_e)\left(X,V_1,V_2,1,-2\pi i,-4\pi^2\right) &= p_e\left(X,\frac{Y_0}{Y_0},\frac{Y_1}{Y_0},\frac{Y_2}{Y_0}\right) \\
        &= Y_0^{-e} p_e(X,\Y)
    \end{align*}
    where the last equality follows from the homogeneity of $p_e$. Since by assumption $p_e \neq 0$, this shows that $r$ is non-trivial, as desired.
\end{proof}

\begin{example}
    It is easy to construct examples where a function has no exponential growth at \emph{some} cusp. For instance, if $p(\Y) = 4\pi^2 Y_0 + Y_2$, then
    \[ p(\J(z)) = 4\pi^2 j(z) + j''(z) \]
    is bounded for $z \to i\infty$ in the standard fundamental domain. However, after the transformation $z \mapsto -\frac{1}{z}$ one gets
    \[ 4\pi^2 j(z) + z^4 j''(z) + 2z^3 j'(z) \]
    which has order $(-1,-4)$ at $i\infty$, since $z^4 j''(z)$ is the dominant term. Note that here $4 = \deg_W(\Gamma(p))$.

    There are also simple examples where the dominant terms cancel out at all cusps. Take $h(\Y) = Y_1^2 - Y_0 Y_2$. Then
    \[ h(\J(\gamma z)) = (cz+d)^4 j'(z)^2 - (cz+d)^4 j(z)j''(z) - 2c(cz+d)^3 j(z) j'(z) \]
    has order $(-2,-3)$ for $c \neq 0$, because $j(z)j''(z) \sim (j'(z))^2 \sim -4\pi^2 q^{-2}$, and it has order $(-2,0)$ for $c = 0$. On the other hand, $\deg_W(\Gamma(h)) = 4$. \Cref{prop:generic order at infinity} guarantees that even though these cancellations may occur at all cusps, some term of maximal exponential growth is not cancelled, at least generically.
\end{example}

Finally we compute the order of $p(z,\J(z))$ at the points $\tau$ for which $j'(\tau) = 0$, namely the orbits of $\rho$ and $i$. Here, the order is at least the maximum $\nu$ such that $T^\nu$ divides respectively $p(X,T^3 Y_0, T^2 Y_1,TY_2)$ (for the conjugates of $\rho$) and $p(X,T^2 Y_0 + 1728,T Y_1,Y_2)$ (for the conjugates of $i$). This estimate is not sharp when $p$ depends on $Y_2$, as we show in an example below, so we restrict to $p \in \C[X,Y_0,Y_1]$.

\begin{proposition}\label{prop:generic order at ramified tau}
    Let $p \in \C[X,Y_0,Y_1]$, $\tau \in \h$, $u = j(\tau)$, and let $\mu$ be the order of $j(z) - u$ at $z = \tau$. Then for all $\gamma$ in a Zariski open dense subset of $\SL_2(\Z)$, the order of $p(z,\J(z))$ at $z = \gamma\tau$ is the highest power of $T$ dividing $p(X,T^\mu Y_0 + u,T^{\mu - 1} Y_1)$.
\end{proposition}
\begin{proof}
    The proof is very similar to that of \Cref{prop:generic order at infinity}. Write
    \[ p(X,T^\mu Y_0 + u,T^{\mu-1}Y_1) = \sum_{k=\nu}^m T^k p_k(X,\Y) \]
    where each $p_k$ satisfies the homogeneity condition
    \[ p_k(X,T^\mu Y_0,T^{\mu-1} Y_1) = T^k p_k(X,Y_0,Y_1) \]
    and $p_\nu \neq 0$. Since
    \[ \Gamma(p)(X,T^\mu Y_0 + u,T^{\mu-1}Y_1) = \Gamma(p(X,T^\mu Y_0 + u,T^{\mu-1}Y_1)), \]
    we have
    \[ \Gamma(p)(T^\mu Y_0 + u,T^{\mu-1}Y_1) = \sum_{k=\nu}^m T^k \Gamma(p_k)(X,\Y) \]
    and each $\Gamma(p_k)$ satisfies the above homogeneity condition.

    Fix $\alpha_0 \in \C$ such that $j(z) - u \sim \alpha_0(z - \tau)^\mu$, and so also $j'(z) \sim \mu\alpha_0(z - \tau)^{\mu-1}$, for $z \to \tau$, where by assumption $\mu\alpha_0 \neq 0$. For $\gamma = \begin{pmatrix} a & b \\ c & d \end{pmatrix}\in \SL_2(\Z)$ we have
    \[ p^\gamma(z,\J(z)) = (c\tau+d)^n r(\gamma\tau,c\tau+d,c)(z - \tau)^\nu + O((z-\tau)^{\nu-1}) \]
    as $z \to \tau$, where
    \[ r(Z,W,C) \coloneqq \Gamma(p_\nu)(Z,W,C,\alpha_0,\mu\alpha_0). \]
    It follows that the order of $p^\gamma(z,\J(z))$ at $z = \tau$ is $\nu$ as long as $r$ does not vanish. Since $z \mapsto \gamma z$ is a diffeomorphism, this coincides with the order of $p(z,\J(z))$ at $z = \gamma\tau$.

    To conclude, since the map $\gamma \mapsto (\gamma \tau,c\tau+d,c)$ is injective on $\SL_2(\C)$, hence dominant on $\C^3$, we only need to show that $r$ is non-trivial. Pick $\beta_0$ such that $\beta_0^{2\mu} = \alpha_0^{-1}$ and $U_0$ such that $U_0^{2\mu} = Y_0$. By \Cref{lem:specialise generic transform} at $\alpha^2 = \mu\alpha_0$, $\beta = 0$, $U_1^2 = (\beta_0 U_0)^{-2(\mu-1)}Y_1$, there are $V_1,V_2 \in \C(U_1,U_2)$ such that
    \begin{align*}
        r(X,V_1,V_2) = \Gamma(p_\nu)(X,V_1,V_2,\alpha_0,\mu\alpha_0) &= p_k\left(X,\alpha_0,\frac{Y_1}{(\beta_0 U_0)^{2(\mu-1)}}\right) \\
        &= p_\nu\left(X,\frac{U_0^{2\mu}}{(\beta_0 U_0)^{2\mu}},\frac{Y_1}{(\beta_0 U_0)^{2(\mu-1)}}\right) \\
        &= (\beta_0 U_0)^{-2k} p_\nu(X,Y_0,Y_1)
    \end{align*}
    where the last equality is implied by the homogeneity condition. It follows that $r$ is non-trivial, as desired.
\end{proof}

\begin{example}\label{example: no zero estimate for j''}
    The above method fails when $p$ depends on $Y_2$. Indeed, to compute say the order of $p(\J(\gamma z))$ at $\rho$, we would look at the maximum power of $T$ dividing $p^\gamma(X,T^3 Y_0,T^2 Y_1,TY_2)$. However, if $p$ contains $Y_2$, then
    \[ \Gamma(p)(X,T^3 Y_0,T^2 Y_1,TY_2) \neq \Gamma(p(X,T^3 Y_0,T^2 Y_1,T Y_2)), \]
    breaking the very first steps of the argument.

    The order can indeed be higher than expected at all the conjugates of $\rho$ or $i$. For instance, for the polynomial
    \[ p(\Y) = Y_0Y_2 - \frac{2}{3}Y_1^2, \]
    the maximum power of $T$ dividing $p(T^3 Y_0,T^2 Y_1,T Y_2)$ is $4$, however the function $p(\J(z))$ has order at least $5$ at all conjugates of $\rho$: if $j(z) \sim \alpha(z - \gamma\rho)^3$ for $z \to \gamma\rho$, then
    \[ p(\J(z)) = \alpha^2(z - \gamma\rho)^4\left(6 - \frac{2}{3}3^2 + O(z - \gamma\rho)\right) = O((z - \gamma\rho)^5). \]
\end{example}

\begin{corollary}\label{cor:special Liouville}
    For all $p \in \C[Y_0,Y_1]$, $h \in \C[Y_0]$, and $\ell \geq \deg_{Y_1}(p)$, if the function
    \[ P(z):=\frac{p(j(z),j'(z))}{h(j(z)) j'(z)^\ell} \]
    is non-constant, then it has a pole in $\h$ or it has exponential growth in some fundamental domains.
\end{corollary}
\begin{proof}
    Suppose that $P(z)$ has no poles in $\h$ and has no exponential growth in any fundamental domain. Without loss of generality, we may assume that $p$ and $h$ are coprime.

    Let $\beta$ be a root of $h(Y_0)$ such that $\beta \notin \{0,1728\}$, and let $\tau$ be such that $j(\tau) = \beta$. In particular, $Y_0 - \beta$ does not divide $p$. By \Cref{prop:generic order at unramified tau}, $p(\gamma\tau,j(\gamma\tau)) \neq 0$ for all $\gamma \in \SL_2(\Z)$ except for some proper Zariski closed set, and so $P$ has a pole at $\gamma\tau$, a contradiction. Therefore, $h(Y_0)$ is of the form $\alpha Y_0^s (Y_0 - 1728)^t$ for some $\alpha \in \C$.

    Now $h(j(z))j'(z)^\ell$ has order $3s + 2\ell$ at all the conjugates of $\rho$ and $2t + \ell$ at all the conjugates of $i$, thus $p(j(z),j'(z))$ must have at least the same order at those points. On writing $p = \sum_{\ell'} Y_1^{\ell'} p_{\ell'}(Y_0)$, \Cref{prop:generic order at ramified tau} implies that each $p_{\ell'}(Y_0)$ is divisible by $Y_0^{s + \lceil\frac{2(\ell - \ell')}{3}\rceil}$ and by $(Y_0 - 1728)^{t + \lceil\frac{\ell - \ell'}{2}\rceil}$. In particular, whenever $p_{\ell'} \neq 0$ we have
    \[ \deg(p) \geq \deg(p_{\ell'}) \geq s + t + \left\lceil\frac{2(\ell - \ell')}{3}\right\rceil + \left\lceil\frac{\ell - \ell'}{2}\right\rceil + \ell' \geq s + t + \ell' + \frac{7}{6}(\ell - \ell'), \]
    with strict inequality if $p_{\ell'}$ is not a constant multiple of $Y_0^{s + \lceil\frac{2(\ell - \ell')}{3}\rceil}(Y_0 - 1728)^{t + \lceil\frac{\ell - \ell'}{2}\rceil}$.

    On the other hand, by \Cref{prop:generic order at infinity}, on a Zariski open dense set of fundamental domains, the denominator has order at most $(-(s + t + \ell),0)$ at the cusp, while the numerator has order at most $(-\deg(p),0)$, thus whenever $p_{\ell'} \neq 0$ we also have
    \[ s + t + \ell' + \frac{7}{6}(\ell - \ell') \leq \deg(p) \leq s + t + \ell. \]
    Since $\ell' \leq \deg_{Y_1}(p) \leq \ell$ and $\frac{7}{6}>1$, it follows at once that $\ell = \ell'$ and that $p = p_{\ell}$ is a constant multiple of $Y_0^s (Y_0 - 1728)^t Y_1^\ell$, and so that $P(z)$ is constant.
\end{proof}

\section{The Main Result}
\label{sec:eq(z,j,j',j'')}

\subsection{\texorpdfstring{$\J$}{j}-Homogeneous equations}
\label{subsec:homeq(j,j'j,'')}

Before tackling the general case of Theorem \ref{thm: intro - main thm}, we look at equations of the form $F(\J(z))=0$ where $F \in \C[\Y]$ satisfies the homogeneity condition below.

\begin{definition}\label{def:j-homogeneous}
    The \textit{$\J$-degree} of $F\in\C[X,\Y]$ is the degree of $F(X,Y_0,T^2 Y_1,T^4 Y_2)$ in $T$, which we denote $\deg_\J(F)$.
    We say that $F$ is \textit{$\J$-homogeneous} if $F(X,Y_0,T^2 Y_1,T^4 Y_2)$ is homogeneous in the variable $T$.
\end{definition}

One of the easiest examples of a $\J$-homogeneous polynomial is $F = Y_2$, which has $\J$-degree 4. For the sake of exposition, we first sketch the proof of Zariski density for this $F$, namely for the equation $j''(z) = 0$.

First, we observe that for $\gamma=\begin{pmatrix}a & b\\ c & d\end{pmatrix}\in\mathrm{SL}_2(\Z)$,
\begin{equation}\label{eq:j''(gamma z)}
    j''(\gamma z) = j''(z) c^4\left(\left(z + \frac{d}{c}\right)^4 + 2\left(z + \frac{d}{c}\right)^3 \frac{j'(z)}{j''(z)}\right) = j''(z) c^4 h(z + \sfrac{d}{c},z)
\end{equation}
where $h(X,W) = X^4 + 2\frac{j'(W)}{j''(W)}X^3$. If we can find $\tau \in \h$ such that $j''(\tau) = 0 \neq j'(\tau)$, we are done by \Cref{cor: pole at cusp Zariski density}, so suppose by contradiction that this does not happen. In this case, $\sfrac{j'(z)}{j''(z)}$ is bounded on the standard fundamental domain $\F$: by construction, it cannot have poles in $\h$, and by looking at the $q$-expansions, it is also bounded for $\im(z) \to +\infty$. This also implies that $h(\tau + \sfrac{d}{c},\tau) \neq 0$ for any $\sfrac{d}{c} \in \Q$, $\tau \in \h$ except possibly when $j'(\tau) = 0$.

Second, under the above assumptions, we shall verify that $h(\tau + r,\tau) \neq 0$ for all $r \in \R$, $\tau \in \h$ (\Cref{claim:h has no zeroes}), and in turn deduce that
\[ \left|h\left(z + \frac{d}{c},z\right)\right| \geq \varepsilon\left|z + \frac{d}{c}\right|^4 \]
for all $z \in \F$, for some $\varepsilon > 0$ (\Cref{claim:h grows polynomially}).

In particular, for all $z \in \F$ we have
\[ \left|\frac{j'(\gamma z)}{j''(\gamma z)}\right| \leq \left|\frac{j'(z)}{\varepsilon j''(z)}\right| |cz + d|^{-2} = \left|\frac{j'(z)}{\varepsilon j''(z)}\right|\frac{\im(\gamma z)}{\im(z)} \leq \frac{2M\im(\gamma z)}{\sqrt{3}\varepsilon} \]
where $M$ is a bound for $\left|\frac{j'(z)}{j''(z)}\right|$ on $\F$, and so $\frac{j'(z)}{j''(z)} \to 0$ as $\im(z) \to 0$. By the Schwarz Reflection Principle, $\frac{j'(z)}{j''(z)}$ extends to a holomorphic function on $\C$ that vanishes for all $z \in \R$, and is thus constantly 0, a contradiction.

For a general $\J$-homogeneous $F$, we just need to find an appropriate generalisation of \Cref{eq:j''(gamma z)} and fill the details in the above sketch.

\begin{lemma}\label{lem:J-homogeneous transform}
    Let $F \in \C[\Y]$ be $\J$-homogeneous. Then there are polynomials $p_k \in \C[\Y]$ and $h \in \C[X,\Y]$ such that for all $\gamma = \begin{pmatrix} a&b \\ c&d \end{pmatrix} \in \SL_2(\Z)$ with $c \neq 0$ we have
    \[ F(\J(\gamma z)) = F(\J(z)) c^N \left(\sum_{k = k_0}^{N} \frac{p_k(\J(z))}{F(\J(z))} \left(z + \frac{d}{c}\right)^k\right) = F(\J(z)) c^N h(z + \sfrac{d}{c},\J(z)), \]
    where $N = \Jdeg(F)$, $p_N = F$, and $0 \neq p_{k_0} \in Y_1^\ell \cdot \C[Y_0]$ with $2\ell \leq k_0$.
\end{lemma}
\begin{proof}
    Let $F$ be as in the hypothesis. By $\J$-homogeneity, we have that
    \begin{align*}
         \Gamma(F)(Z,W,C,\Y) &= F(Z,Y_0,W^2 Y_1,W^4 Y_2 + 2CW^3 Y_1) \\
         &= C^N F\left(Z,Y_0,\frac{W^2}{C^2}Y_1,\frac{W^4}{C^4} Y_2 + 2\frac{W^3}{C^3} Y_1\right).
    \end{align*}
    Therefore, we can write $\Gamma(F)$ as
    \[ \Gamma(F) = C^N \sum_{k = 0}^{N} p_k(\Y)\frac{W^k}{C^k} \]
    where $N = \Jdeg(F)$. Let $k_0$ be the least integer such that $p_{k_0} \neq 0$.

    By a further application of $\J$-homogeneity, we also have
    \[ \Gamma(F) = W^N Y_1^{\frac{N}{2}} F\left(Z,Y_0,1,\frac{Y_2}{Y_1^2} + \frac{2C}{WY_1}\right). \]
    It follows at once that the terms of maximum degree in $W$, which make up $\frac{W^N}{C^N}p_N$, are found by discarding $\frac{2C}{WY_1}$, and in particular we discover that $p_N(\Y) = F(\Y)$. Similarly, the terms of lowest degree in $\frac{W}{C}$ are obtained by specialising at $Y_2 = 0$ and taking the least power of $Y_1$, and so $p_{k_0} \in Y_1^\ell \cdot \C[Y_0]$. Moreover, if $t$ is the degree of $F$ in $Y_2$, we have that $k_0 = N - t$, $\ell = \frac{N}{2} - t$, thus $2\ell \leq k_0$.
\end{proof}

\begin{theorem}\label{thm:homogeneous}
    For any irreducible $\J$-homogeneous polynomial $F(\Y) \notin \C[X,Y_0,Y_1]$, the equation $F(\J(z)) = 0$ has a Zariski dense set of solutions.
\end{theorem}
\begin{proof}
    Let $F$ be as in the hypothesis and fix the polynomials $p_k$, $h$ as in the conclusion of \Cref{lem:J-homogeneous transform}.

    If some $\frac{p_k(\J(z))}{F(\J(z))}$ has a pole in $\h$ or exponential growth in some fundamental domain, we are done by \Cref{cor: pole at cusp Zariski density}. Therefore, we shall assume that this is not the case. In particular, we assume that $h(z + \sfrac{d}{c},\J(z))$ has no pole in $\h$ for any $\sfrac{d}{c}$.

    With this additional assumption, if $F(\J(\tau)) = 0$, then $F(\J(\gamma\tau)) = 0$ for all $\gamma \in \SL_2(\Z)$. Since $F$ is not divisible by $Y_0 - j(\tau)$, \Cref{prop:generic order at unramified tau} implies that $j'(\tau) = 0$. In particular, for any $\sfrac{d}{c}$, $h(\tau + \sfrac{d}{c},\J(\tau)) = 0$ implies that $j'(\tau) = 0$.

    \begin{claim}\label{claim:h has no zeroes}
        For every $u \in \R$, the function $h(z + u,\J(z))$ has no zero in $\h$.
    \end{claim}
    \begin{claimproof}
        Suppose $h(\tau + u,\J(\tau)) = 0$ for some $(\tau,u) \in \h \times \R$. Since $h(Z,\Y)$ is monic in $Z$, the analytic map $(z,u) \mapsto (h(z + u,\J(z)),u)$ has finite fibres, in particular of dimension zero. Then the Open Mapping Theorem implies that the image of any ball around $(\tau,u)$ contains an open neighbourhood of $(0,u)$. In particular, for every rational number $r$ arbitrarily close to $u$, there is $\tau_r$ close to $\tau$ such that $h(\tau_r + r,\J(\tau_r)) = 0$.

        Since $h$ is monic in $Z$, the polynomial $h(\tau + Z,\J(\tau))$ is not identically zero, thus for $r$ sufficiently close to $u$ we have $h(\tau + r,\J(\tau)) \neq 0$, and so $\tau_r \neq \tau$. Therefore, the $\tau_r$'s can be chosen to accumulate at $\tau$ for $r \to u$. However, our assumptions imply that $j'(\tau_r) = 0$, so the $\tau_r$'s lie in a closed discrete subset of $\h$ (namely, the orbits of $\rho$ and $i$), a contradiction.
    \end{claimproof}

    \begin{claim}\label{claim:h grows polynomially}
        There is $\varepsilon > 0$ such that
        \[ |h(z + u,\J(z))| \geq \varepsilon |z + u|^N \]
        for all $z \in \F$ and $u \in \R$.
    \end{claim}
    \begin{claimproof}
        Our current assumptions imply that each $\frac{p_k(\J(z))}{F(\J(z))}$ has neither a pole in $\overline{\F}$ nor exponential growth in $\F$. Since these functions are in $\mathcal{P}$ rather than $\mathcal{P}(w)$, their order at $i\infty$ is of the form $(e,0)$, thus at least $(0,0)$, and so they are bounded in $\overline{\F}$, say by $M > 0$. In particular, since $h$ is monic in $Z$, we have $|h(z+u,\J(z))| > \frac{1}{2} |z + u|^N$ as soon as $|z + u| > 2N\!M$. On the other hand, $|z + u| \leq 2N\!M$ defines a compact subset $K$ of $\overline{\F} \times \R$, and so the function
        \[ \left|\frac{h(z + u,\J(z))}{(z + u)^N}\right| \]
        attains some minimum $\varepsilon' > 0$ on $K$, since it does not vanish by \Cref{claim:h has no zeroes}. The conclusion follows on taking $\varepsilon = \min\{\varepsilon',\frac{1}{2}\}$.
    \end{claimproof}

    Therefore, for $z \in \F$, we find that
    \[ |F(\J(\gamma z))| \geq \varepsilon |F(\J(z))| c^N \left|z + \frac{d}{c}\right|^N = \varepsilon |F(\J(z))| |cz + d|^N, \]
    and in particular for some $M > 0$ independent of $z$ we have
    \[ \left|\frac{p_{k_0}(\J(\gamma z))}{F(\J(\gamma z))}\right| = \frac{|p_{k_0}(\J(z))| |cz+d|^{2\ell}}{|F(\J(\gamma z))|} \leq \frac{|p_{k_0}(\J(z))|}{\varepsilon|F(\J(z))|} |cz+d|^{2\ell - N} \leq M \im(\gamma z)^{\ell - \frac{N}{2}}, \]
    where we have used that $2\ell < N$, and that for $z \in \F$, $\frac{p_{k_0}(\J(z))}{F(\J(z))}$ is bounded while also
    \[ |cz+d|^2 = \frac{\im(z)}{\im(\gamma z)} \geq \frac{\sqrt{3}}{2\im(\gamma z)} > 0. \]
    Therefore, $\frac{p_{k_0}(\J(z))}{F(\J(z))} \to 0$ as $\im(z) \to 0$, thus by Schwarz's reflection principle, $\frac{p_{k_0}(\J(z))}{F(\J(z))}$ has a holomorphic extension to $\C$ that is constantly zero on $\R$, thus constantly zero on $\C$, hence $p_{k_0} = 0$, a contradiction.
\end{proof}

\begin{remark}\label{rem:j(S) is almost perfect}
    Let $S = \{ z \in \h : F(\J(z)) = 0 \}$ for some $F$ as in \Cref{thm:homogeneous}. As observed in the proof, every $\tau \in S$ that is not in the $\SL_2(\Z)$-orbit of $\rho$ or $i$ must also be a pole of some coefficient $\frac{p_k(\J(z))}{F(\J(z))}$ appearing in \Cref{lem:J-homogeneous transform}. Combining this information with \Cref{prop: periodic finite pole}, it follows that $j(S)$ is infinite and every point of $j(S) \setminus \{0,1728\}$ is an accumulation point of $j(S)$.
\end{remark}

\subsection{Proof of Theorem \ref{thm: intro - main thm}}
\label{subsec:proofmainthm}

For the rest of the section, fix some $F\in\C[X,\Y]$ and let $p_k \in \C[Z,C,\Y]$ be polynomials such that
\[ \Gamma(F)(Z,W,C,\Y) = \sum_{k = 0}^{N} p_k(Z,C,\Y)W^k, \]
where $N = \deg_W(\Gamma(F))$.
Given a polynomial $\alpha\in \C[X,\Y]$, we will use the notation\footnote{We remark that for us $0\in\N$.}
\begin{equation*}
    \alpha^\N :=\{\alpha^s : s\in\N\}.
\end{equation*}
Given different polynomials $\alpha_1,\ldots,\alpha_\ell\in\C[X,\Y]$, we use $\alpha_1^\N\cdots\alpha_\ell^\N$ to denote the set of all products between elements of different $\alpha_t^\N$.

\begin{proposition}\label{prop:pN is homogeneous}
    The polynomial $p_N$ is the sum of the terms of maximum $\J$-degree in $F(Z,Y_0,W^2 Y_1,W^4 Y_2)$. In particular, $N = \Jdeg(F)$, $p_N$ is $\J$-homogeneous, and $p_N$ does not depend on $C$.
\end{proposition}
\begin{proof}
   { Let $X^\alpha Y_0^{\beta_0}Y_1^{\beta_1}Y_2^{\beta_2}$ denote a monomial, so $\alpha,\beta_0,\beta_1,\beta_2\in\N$.
    Observe that
    \[\Gamma(X^\alpha Y_0^{\beta_0}Y_1^{\beta_1}Y_2^{\beta_2}) = Z^\alpha Y_0^{\beta_0}Y_1^{\beta_1}W^{2\beta_1 + 3\beta_2}(WY_2 + 2CY_1)^{\beta_2}.\]
    Hence
    \[\deg_W(\Gamma(X^\alpha Y_0^{\beta_0}Y_1^{\beta_1}Y_2^{\beta_2})) = 2\beta_1 + 4\beta_2\]
    and the term accompanying $W^{2\beta_1+4\beta_2}$ is $Z^\alpha Y_0^{\beta_0}Y_1^{\beta_1}Y_2^{\beta_2}$.

    Now write $\boldsymbol{\beta}=(\beta_0,\beta_1,\beta_2)$ and
    \[F(X,\Y) = \sum_{(\alpha,\boldsymbol{\beta})\in\N^4}c_{\alpha,\boldsymbol{\beta}}X^\alpha\Y^{\boldsymbol{\beta}},\]
    where $\Y^{\boldsymbol{\beta}} = Y_0^{\beta_0}Y_1^{\beta_1}Y_2^{\beta_2}$ and $c_{\alpha,\boldsymbol{\beta}}\in\C$.
    Then, since $\Gamma$ is a homomorphism by \eqref{Gamma-homomorphism},
    \[\Gamma(F) = \sum_{(\alpha,\boldsymbol{\beta})\in\N^4}c_{\alpha,\boldsymbol{\beta}}\Gamma\left(X^\alpha\Y^{\boldsymbol{\beta}}\right).\]
    Since $N=\deg_W(\Gamma(F))$, then
    \[p_N = \sum_{(\alpha,\boldsymbol{\beta})\in\N^4: 2\beta_1+4\beta_2=N}c_{\alpha,\boldsymbol{\beta}}Z^\alpha Y_0^{\beta_0}Y_1^{\beta_1}Y_2^{\beta_2}.\]
    From this we see that $p_N$ does not depend on $C$ and that $p_N$ is $\J$-homogeneous.
    This expression also gives us that $p_N$ is the coefficient of $W^N$ in $F(Z,Y_0,W^2 Y_1,W^4 Y_2)$, and so $N=\deg_W(F(Z,Y_0,W^2 Y_1,W^4 Y_2)) = \deg_\J(F)$.}
\end{proof}

\begin{definition}
    The \textit{$\J$-order} of $F$, denoted by $\Jord(F)$, is the maximum power of $T$ dividing $F(X,Y_0,T^2 Y_1,T^3 Y_2)$.
\end{definition}

\begin{proposition}\label{prop:pk has no Y2}
    Let $k_0$ be minimum such that $p_{k_0} \neq 0$. Then $p_{k_0}$ is the sum of the terms of minimum degree in $W$ of $F(Z,Y_0,W^2 Y_1, 2CW^3 Y_1)$. In particular, $k_0 = \Jord(F) \geq 2\deg_{Y_1}(p_{k_0})$ and $p_{k_0}$ does not depend on $Y_2$.
\end{proposition}
\begin{proof}
    We proceed as in the proof of \Cref{prop:pN is homogeneous}.
    From
    \[\Gamma(X^\alpha Y_0^{\beta_0}Y_1^{\beta_1}Y_2^{\beta_2}) = Z^\alpha Y_0^{\beta_0}Y_1^{\beta_1}W^{2\beta_1 + 3\beta_2}(WY_2 + 2CY_1)^{\beta_2}\]
    we see that the smallest power of $W$ appearing in this expression is $2\beta_1+3\beta_2$, and it is accompanied by
    \[2^{\beta_2}C^{\beta_2}Z^\alpha Y_0^{\beta_0}Y_1^{\beta_1+\beta_2},\]
    which does not depend on $Y_2$.
    Hence
    \[p_{k_0}=\sum_{(\alpha,\boldsymbol{\beta})\in\N^4: 2\beta_1+3\beta_2=k_0}c_{\alpha,\boldsymbol{\beta}}2^{\beta_2}C^{\beta_2}Z^\alpha Y_0^{\beta_0}Y_1^{\beta_1+\beta_2},\]
    which also shows that $p_{k_0}$ is the coefficient of $W^{k_0}$ in $F(Z,Y_0,W^2 Y_1, 2CW^3 Y_1)$.
    Using the change of variables $C = \frac{Y_2}{2Y_1}$, we conclude that $k_0=\deg_\J(F)$, and since $k_0 = 2\beta_1+3\beta_2\geq 2(\beta_1+\beta_2)$ we have $k_0 \geq 2\deg_{Y_1}(p_{k_0})$, concluding the proof.
\end{proof}

\begin{corollary}\label{cor:pN beats pk on Y1}
    If $F \notin Y_1^\N \C[X,Y_0]$, then $\Jdeg(F) > 2\deg_{Y_1}(p_{\Jord(F)})$.
\end{corollary}
\begin{proof}
    One can immediately verify that
    \begin{align*}
      \Jdeg(F) &= \deg_T(F(X,Y_0,T^2 Y_1,T^4 Y_2))\\
      &\geq \deg_T(F(X,Y_0,T^2 Y_1,T^3 Y_2)) \\
      &\geq \ord_{T=0}(F(X,Y_0,T^2 Y_1,T^3 Y_2)) \\
      &= \Jord(F) \\
      &\geq 2\deg_{Y_1}(p_{\ord_\J(F)}),
    \end{align*}
    where $\ord_{T=0}(P)$ is the maximum power of $T$ dividing $P$.

    If the second inequality is an equality, then $F(X,Y_0,T^2 Y_1,T^3 Y_2) = T^m F(X,Y_0,Y_1,Y_2)$ where $m = \Jord(F)$. In particular,
    \[ \deg_T(F(X,Y_0,T^2 Y_1, T^4 Y_2)) = \deg_T(T^m F(X,Y_0,Y_1,T Y_2)) = m + \deg_{Y_2}(F). \]
    If the first inequality is also an equality, then $\deg_{Y_2}(F) = 0$, and moreover $F$ is homogeneous in $Y_1$ of degree $\frac{m}{2}$, thus $F \in Y_1^\N \C[X,Y_0]$.
\end{proof}

We can now prove the main result of this paper, Theorem \ref{thm: intro - main thm}, the statement of which is recalled below for the convenience of the reader.

\begin{customthm}{\ref{thm: intro - main thm}}
     For any polynomial $F(X,\Y) \in \C[X, \Y]\setminus \C[X]$ which is coprime to $Y_0(Y_0-1728)Y_1$, the equation $F(z,\J(z))=0$ has a Zariski dense set of solutions.
 \end{customthm}
\begin{proof}
    It suffices to prove the conclusion for $F$ irreducible, not in $\C[X]$, and not a constant multiple of $Y_0$, $Y_0 - 1728$, or $Y_1$. Let $n = \deg_X(F)$ and write
    \begin{align*}
       F^\Gamma &= (CX+D)^n \Gamma(F)\left(\frac{AX+B}{CX+D},CX+D,C,\Y\right) \\
                &= (CX+D)^n \sum_{k=0}^N p_k\left(\frac{AX+B}{CX+D},C,\Y\right) (CX+D)^k \\
                &= \sum_{k=0}^{n+N} h_k(A,B,C,D,\Y)X^k.
    \end{align*}
    We recall that $F^\gamma = F^\Gamma(a,b,c,d,X,\Y)$ for any $\gamma=  \begin{pmatrix} a & b \\ c & d \end{pmatrix}\in \SL_2(\Z)$. If $F \in \C[Y_0]$, the conclusion follows by \Cref{cor:j(z)-u Zariski dense}, so we may assume that this is not the case, and in particular that $n + N > 0$.

    We observe immediately that $h_{n+N} = C^n p_N(\sfrac{A}{C},\Y)$ using \Cref{prop:pN is homogeneous}. Moreover, by \Cref{prop:pk has no Y2},
    \[ h_0 = D^n \sum_{k=0}^N p_k(\sfrac{B}{D},C,\Y) D^k = D^n \sum_{k=\Jord(F)}^N p_k(\sfrac{B}{D},C,\Y) D^k. \]
    We claim that for $\gamma= \begin{pmatrix} a & b \\ c & d \end{pmatrix}$ in some Zariski open dense subset of $\SL_2(\Z)$, we have that $h_{n+N}(a,b,c,d,\Y)$ has a factor $r\in\C[\Y]$ such that the equation $r(\J(z)) = 0$ has a Zariski dense set of solutions, or that $\frac{h_0}{h_{n+N}}(a,b,c,d,\J(z))$ has a pole or exponential growth in some fundamental domain. In particular, for any one of those $\gamma$'s we find that
    \[ F^\gamma(z,\J(z)) = \sum_{k=0}^{n+N} h_k(a,b,c,d,\J(z)) z^k = 0 \]
    has a Zariski dense set of solutions (by \Cref{cor:Zariski density begets Zariski density} in the first case, and \Cref{cor: pole at cusp Zariski density} in the second one), hence so does $F(z,\J(z)) = 0$ (by Proposition \ref{prop:p^gamma}), as desired.

    To prove the claim, we distinguish three cases.

    \medskip
    \noindent\textbf{Suppose that $p_N$ is not in $Y_0^\N (Y_0 - 1728)^\N Y_1^\N \C[Z]$.} Recall that $p_N$ is $\J$-homogeneous by \Cref{prop:pN is homogeneous}. Since $h_{n+N} = C^n p_N(\sfrac{A}{C},\Y)$, if $p_N$ depends on $Y_2$, so does $h_{n+N}(\sfrac{a}{c},\Y)$ for all but finitely many values of $\sfrac{a}{c}$.

    If $p_N$ does not depend on $Y_2$, then $p_N = Y_1^N h(\sfrac{A}{C},Y_0)$ for some polynomial $h\in\C[Z,Y_0]$, and by assumption $h(\sfrac{A}{C},Y_0)$ has at least one root distinct from $0$ and $1728$, seen as a polynomial in $Y_0$ (in an algebraic closure of $\C(\sfrac{A}{C})$), in which case so does $h(\sfrac{a}{c},Y_0)$ except for finitely many values of $\sfrac{a}{c}$.

    In either case, since the map $\gamma \mapsto \sfrac{a}{c}$ from $\SL_2(\C)$ to $\C$ is dominant, we get that for all $\gamma$ except on some proper Zariski closed subset of $\SL_2(\Z)$, the polynomial $h_{n+N}(a,b,c,d,\Y) = c^n p_N(\sfrac{a}{c},\Y)$ has an irreducible factor $r\in\C[\Y]$ such that $r(\J(z)) = 0$ has a Zariski dense set of solutions (by respectively \Cref{thm:homogeneous}, after noticing that the factors of a $\J$-homogeneous polynomial are $\J$-homogeneous, and \Cref{cor:j(z)-u Zariski dense}).

    \medskip
    \noindent\textbf{Suppose that $F$ is in $\C[X,Y_0]$.} By irreducibility of $F$, we have that $F = F(X,Y_0)$ is neither divisible by $Y_0$ nor $Y_0 - 1728$. In particular, $p_N = \Gamma(F) = F(Z,Y_0)$ is also neither divisible by $Y_0$ nor $Y_0 - 1728$, which puts us back in the previous case.

    \medskip
    \noindent\textbf{Suppose that $p_N$ is in $Y_0^\N (Y_0 - 1728)^\N Y_1^\N \C[Z]$ and $F$ is not in $\C[X,Y_0]$.} Write $p_N = r(Z)Y_0^s (Y_0 - 1728)^t Y_1^\ell$. Since $F$ is irreducible, $F$ is also not in $Y_1^\N \C[X,Y_0]$, thus by \Cref{cor:pN beats pk on Y1},
    \[ 2\ell = 2\deg_{Y_1}(p_N) = \Jdeg(p_N) = \Jdeg(F) > 2\deg_{Y_1}(p_{\Jord(F)}). \]
    Since the map $\gamma \mapsto (\sfrac{a}{c},\sfrac{b}{d},c)$ from $\SL_2(\C)$ to $\C^3$ is dominant, for $\gamma$ in some Zariski open dense subset of $\SL_2(\Z)$ we have
    \[ \deg_{Y_1}(p_N(\sfrac{a}{c},c)) = \deg_{Y_1}(p_N) > \deg_{Y_1}(p_{\Jord(F)}) = \deg_{Y_1}(p_{\Jord(F)}(\sfrac{b}{d},c,\Y)), \]
    in which case the function $\frac{p_{\Jord(F)}(\sfrac{b}{d},c,\J(z))}{p_N(\sfrac{a}{c},\J(z))}$ has a pole at some $\tau \in \SL_2(\Z)\rho \cup \SL_2(\Z)i$ or exponential growth in some fundamental domain by \Cref{cor:special Liouville} (recall that by Proposition \ref{prop:pk has no Y2} we know that $p_{\Jord(F)}$ does not depend on $Y_2$). If we fix some $\sfrac{b}{d}$, $c$ as above, and a corresponding pole $\tau$ or a fundamental domain $\eta\F$ with exponential growth, then the function
    \[ \frac{h_0(\sfrac{b}{d},c,\J(z))}{h_{n+N}(\sfrac{a}{c},\J(z))} = \frac{d^n}{c^n r(\sfrac{a}{c})} \sum_{k=\Jord(F)}^N \frac{p_k(\sfrac{b}{d},c,\J(z))}{j(z)^s (j(z) - 1728)^t j'(z)^\ell}d^k, \]
    has neither a pole at $\tau$ nor exponential growth in $\eta\F$ only when $d$ satisfies a non-trivial polynomial equation over $\sfrac{b}{d}$, $c$. Since $\gamma \mapsto (\sfrac{b}{d},c,d)$ is also a dominant map, for $\gamma$ in a Zariski open dense subset of $\SL_2(\Z)$ the above function has a pole at $\tau$ or exponential growth in $\eta\F$, as claimed.
\end{proof}

\subsection{Two more examples}

\begin{example}
    Let us apply \Cref{thm: intro - main thm} to get some information on the zeroes of the function $j'''$. From the differential equation of the $j$-function (\ref{eq:j}) we see that
    \[ j''' = \frac{3}{2}\cdot \frac{(j'')^2}{j'} - \frac{j^{2} -1968j + 2654208}{2j^{2}(j-1728)^{2}}\left(j'\right)^{3}. \]

    By \Cref{thm: intro - main thm}, the equation
    \begin{equation}\label{eq: j'''=0}
        3 j^2(j-1728)^2 (j'')^2 - \left( j^{2} -1968j + 2654208 \right) (j')^4=0
    \end{equation}
    has a Zariski dense set of solutions outside $\SL_2(\Z)\rho \cup \SL_2(\Z)i$ (this actually follows directly from \Cref{thm:homogeneous}, because the equation is $\J$-homogeneous). Then these are also solutions of $j'''(z) = 0$.

    Upon applying the transformation $z\mapsto -\frac{1}{z}$ we get
    \begin{multline}\label{eq: j'''=0 with z}
         z^8  \left( 3 j(z)^2(j(z)-1728)^2 j''(z)^2 - \left( j(z)^{2} -1968j(z) + 2654208 \right) j'(z)^4 \right)  \\
        +  z^7  12j(z)^2(j(z)-1728)^2 j'(z)j''(z) + z^6  12 j(z)^2(j(z)-1728)^2 j'(z)^2 = 0.
    \end{multline}

    We can see that the ratios \[ \frac{12j^2(j-1728)^2 j'j''}{3 j^2(j-1728)^2 (j'')^2 - \left( j^{2} -1968j + 2654208 \right) (j')^4}\] and \[ \frac{12 j^2(j-1728)^2 (j')^2}{3 j^2(j-1728)^2 (j'')^2 - \left( j^{2} -1968j + 2654208 \right) (j')^4}\] are equal to $0$ at $i$ and $\rho$ and do not have exponential growth in any fundamental domain. However, we know that \eqref{eq: j'''=0} has a zero $\tau \notin \SL_2(\Z)\rho \cup \SL_2(\Z)i$. Therefore, the second ratio above has a pole at $\tau$, and so for all large enough $m$ the equation \eqref{eq: j'''=0 with z} has a zero near $\tau+m$. This means that \eqref{eq: j'''=0}, and hence $j'''=0$, has solutions near $-\frac{1}{\tau+m}$ which accumulate at $0$.
\end{example}

\begin{example}
    Given an equation $F(z,\J(z))=0$, our strategy of the proof of \Cref{thm: intro - main thm} is to apply a \textit{generic} $\SL_2(\Z)$-transformation and show that in the function $F(\gamma z, \J(\gamma z))$ the ratio of the coefficients of the lowest and highest powers of $z$ has a pole at some point in $\tau$ or has exponential growth at a cusp. In some cases, e.g.\ when $F(X,\mathbf{Y})$ does not depend on $X$, we can keep things simple and just use the good old transformation $z\mapsto -\frac{1}{z}$. Indeed, in this case it turns out that the coefficient of the lowest power of $z$ does not depend on $j''$ and so we can apply \Cref{cor:special Liouville}. For instance,
    \[ j'' \left( -\frac{1}{z} \right)^2 = z^8 j''(z)^2 + 4z^7 j'(z) j''(z) + 4z^6 j'(z)^2.\]

    We give an example to show that the transformation $-\frac{1}{z}$ and even $a-\frac{1}{z}$ for any integer $a$ does not suffice in general (when $F$ depends on $X$), as the aforementioned coefficient may depend on $j''$. Consider the function $f(z) = 2j'(z)+zj''(z)$. After a $z \mapsto a-\frac{1}{z}$ transformation we get
    \[ f\left( a-\frac{1}{z}\right) = 2z^2 j'(z) +\left( a-\frac{1}{z}\right) (z^4j''(z)+2z^3j'(z)) = az^4 j''(z) + z^3(2aj'(z)-j''(z)). \]
    We see that the coefficient of $z^3$ depends on $j''$ regardless of the value of $a$.
\end{example}

\subsection*{Acknowledgement} We thank the referee for a thorough reading of the paper and for numerous suggestions that helped us improve the presentation.

\bibliographystyle{alphaurl}
\bibliography{mscrefs}

\begin{thebibliography}{AKM23}

\bibitem[AEK21]{AEK2021}
Vahagn Aslanyan, Sebastian Eterovi\'{c}, and Jonathan Kirby.
\newblock Differential existential closedness for the {$j$}-function.
\newblock {\em Proc. Amer. Math. Soc.}, 149(4):1417--1429, 2021.
\newblock \href {https://doi.org/10.1090/proc/15333} {\path{doi:10.1090/proc/15333}}.

\bibitem[AEK23]{AEK2023}
Vahagn Aslanyan, Sebastian Eterovi\'{c}, and Jonathan Kirby.
\newblock A closure operator respecting the modular {$j$}-function.
\newblock {\em Israel J. Math.}, 253(1):321--357, 2023.
\newblock \href {https://doi.org/10.1007/s11856-022-2362-y} {\path{doi:10.1007/s11856-022-2362-y}}.

\bibitem[AK22]{AK2022}
Vahagn Aslanyan and Jonathan Kirby.
\newblock Blurrings of the {$J$}-function.
\newblock {\em Q. J. Math.}, 73(2):461--475, 2022.
\newblock \href {https://doi.org/10.1093/qmath/haab037} {\path{doi:10.1093/qmath/haab037}}.

\bibitem[AKM23]{AKM2023}
Vahagn Aslanyan, Jonathan Kirby, and Vincenzo Mantova.
\newblock A geometric approach to some systems of exponential equations.
\newblock {\em Int. Math. Res. Not. IMRN}, (5):4046--4081, 2023.
\newblock \href {https://doi.org/10.1093/imrn/rnab340} {\path{doi:10.1093/imrn/rnab340}}.

\bibitem[Asl22]{Asl2022}
Vahagn Aslanyan.
\newblock Adequate predimension inequalities in differential fields.
\newblock {\em Ann. Pure Appl. Logic}, 173(1):Paper No. 103030, 42, 2022.
\newblock \href {https://doi.org/10.1016/j.apal.2021.103030} {\path{doi:10.1016/j.apal.2021.103030}}.

\bibitem[BM17]{BM2017}
W.~Dale Brownawell and David~W. Masser.
\newblock Zero estimates with moving targets.
\newblock {\em J. Lond. Math. Soc. (2)}, 95(2):441--454, 2017.
\newblock \href {https://doi.org/10.1112/jlms.12014} {\path{doi:10.1112/jlms.12014}}.

\bibitem[DFT18]{DFT2018}
Paola D'Aquino, Antongiulio Fornasiero, and Giuseppina Terzo.
\newblock Generic solutions of equations with iterated exponentials.
\newblock {\em Trans. Amer. Math. Soc.}, 370(2):1393--1407, 2018.
\newblock \href {https://doi.org/10.1090/tran/7206} {\path{doi:10.1090/tran/7206}}.

\bibitem[EH21]{EH2021}
Sebastian Eterovi\'{c} and Sebasti\'{a}n Herrero.
\newblock Solutions of equations involving the modular {$j$} function.
\newblock {\em Trans. Amer. Math. Soc.}, 374(6):3971--3998, 2021.
\newblock \href {https://doi.org/10.1090/tran/8244} {\path{doi:10.1090/tran/8244}}.

\bibitem[Ete25]{Ete2022}
Sebastian Eterovi\'c.
\newblock Generic solutions of equations involving the modular {$j$} function.
\newblock {\em Math. Ann.}, 391(4):6401--6449, 2025.
\newblock \href {https://doi.org/10.1007/s00208-024-03082-6} {\path{doi:10.1007/s00208-024-03082-6}}.

\bibitem[Gal23]{Gal2023}
Francesco~Paolo Gallinaro.
\newblock Exponential sums equations and tropical geometry.
\newblock {\em Selecta Math. (N.S.)}, 29(4):Paper No. 49, 41, 2023.
\newblock \href {https://doi.org/10.1007/s00029-023-00853-y} {\path{doi:10.1007/s00029-023-00853-y}}.

\bibitem[Gal25]{Gal2021}
Francesco~Paolo Gallinaro.
\newblock Solving systems of equations of raising-to-powers type.
\newblock {\em Israel J. Math.}, 2025.
\newblock \href {https://doi.org/10.1007/s11856-025-2778-2} {\path{doi:10.1007/s11856-025-2778-2}}.

\bibitem[Kir19]{Kir2019}
Jonathan Kirby.
\newblock Blurred complex exponentiation.
\newblock {\em Selecta Math. (N.S.)}, 25(5):Paper No. 72, 15, 2019.
\newblock \href {https://doi.org/10.1007/s00029-019-0517-4} {\path{doi:10.1007/s00029-019-0517-4}}.

\bibitem[Lan87]{lang-ellipticfunctions}
Serge Lang.
\newblock {\em Elliptic functions}, volume 112 of {\em Graduate Texts in Mathematics}.
\newblock Springer-Verlag, New York, second edition, 1987.
\newblock With an appendix by J. Tate.
\newblock \href {https://doi.org/10.1007/978-1-4612-4752-4} {\path{doi:10.1007/978-1-4612-4752-4}}.

\bibitem[Lan99]{lang:complexanalysis}
Serge Lang.
\newblock {\em Complex analysis}, volume 103 of {\em Graduate Texts in Mathematics}.
\newblock Springer-Verlag, New York, fourth edition, 1999.
\newblock \href {https://doi.org/10.1007/978-1-4757-3083-8} {\path{doi:10.1007/978-1-4757-3083-8}}.

\bibitem[Mah69]{mahler}
Kurt Mahler.
\newblock On algebraic differential equations satisfied by automorphic functions.
\newblock {\em J. Austral. Math. Soc.}, 10:445--450, 1969.
\newblock \href {https://doi.org/10.1017/S1446788700007709} {\path{doi:10.1017/S1446788700007709}}.

\bibitem[Zil02]{Zil2002}
Boris Zilber.
\newblock Exponential sums equations and the {S}chanuel conjecture.
\newblock {\em J. London Math. Soc. (2)}, 65(1):27--44, 2002.
\newblock \href {https://doi.org/10.1112/S0024610701002861} {\path{doi:10.1112/S0024610701002861}}.

\bibitem[Zil05]{Zil2005}
Boris Zilber.
\newblock Pseudo-exponentiation on algebraically closed fields of characteristic zero.
\newblock {\em Ann. Pure Appl. Logic}, 132(1):67--95, 2005.
\newblock \href {https://doi.org/10.1016/j.apal.2004.07.001} {\path{doi:10.1016/j.apal.2004.07.001}}.

\bibitem[Zil15]{Zil2015}
Boris Zilber.
\newblock The theory of exponential sums, 2015.
\newblock \href {https://arxiv.org/abs/1501.03297} {\path{arXiv:1501.03297}}.

\end{thebibliography}

\end{document}